\documentclass[12pt]{article}
\usepackage{amsmath,amsfonts,amssymb,amsthm}
\input amssym.def
\begin{document}
\def \Z{\mathbb Z}
\def \C{\mathbb C}
\def \R{\mathbb R}
\def \Q{\mathbb Q}
\def \N{\mathbb N}

\def \A{{\mathcal{A}}}
\def \D{{\mathcal{D}}}
\def \E{{\mathcal{E}}}
\def \E{{\mathcal{E}}}
\def \H{\mathcal{H}}
\def \S{{\mathcal{S}}}
\def \U{{\mathcal{U}}}

\def \wt{{\rm wt}}
\def \tr{{\rm tr}}
\def \add{{\rm add}}
\def \span{{\rm span}}
\def \Res{{\rm Res}}
\def \Der{{\rm Der}}
\def \End{{\rm End}}
\def \Ind {{\rm Ind}}
\def \Irr {{\rm Irr}}
\def \Aut{{\rm Aut}}
\def \GL{{\rm GL}}
\def \Hom{{\rm Hom}}
\def \mod{{\rm mod}}
\def \ann{{\rm Ann}}
\def \ad{{\rm ad}}
\def \rank{{\rm rank}\;}
\def \<{\langle}
\def \>{\rangle}

\def \g{{\frak{g}}}
\def \h{{\hbar}}
\def \k{{\frak{k}}}
\def \sl{{\frak{sl}}}
\def \gl{{\frak{gl}}}

\def \be{\begin{equation}\label}
\def \ee{\end{equation}}
\def \bex{\begin{example}\label}
\def \eex{\end{example}}
\def \bl{\begin{lem}\label}
\def \el{\end{lem}}
\def \bt{\begin{thm}\label}
\def \et{\end{thm}}
\def \bp{\begin{prop}\label}
\def \ep{\end{prop}}
\def \br{\begin{rem}\label}
\def \er{\end{rem}}
\def \bc{\begin{coro}\label}
\def \ec{\end{coro}}
\def \bd{\begin{de}\label}
\def \ed{\end{de}}

\newcommand{\q}{\bf q}
\newcommand{\m}{\bf m}
\newcommand{\n}{\bf n}
\newcommand{\nno}{\nonumber}
\newcommand{\nord}{\mbox{\scriptsize ${\circ\atop\circ}$}}
\newtheorem{thm}{Theorem}[section]
\newtheorem{prop}[thm]{Proposition}
\newtheorem{coro}[thm]{Corollary}
\newtheorem{conj}[thm]{Conjecture}
\newtheorem{example}[thm]{Example}
\newtheorem{lem}[thm]{Lemma}
\newtheorem{rem}[thm]{Remark}
\newtheorem{de}[thm]{Definition}
\newtheorem{hy}[thm]{Hypothesis}
\makeatletter \@addtoreset{equation}{section}
\def\theequation{\thesection.\arabic{equation}}
\makeatother \makeatletter

\begin{center}
{\Large \bf $\phi$-coordinated modules for quantum vertex algebras and associative algebras}
\end{center}

\begin{center}
{Haisheng Li\\
Department of Mathematical Sciences\\
 Rutgers University, Camden, NJ
08102}
\end{center}

\begin{abstract}
We study $\N$-graded $\phi$-coordinated modules for a general quantum vertex algebra
$V$ of a certain type in terms of an associative algebra $\widetilde{A}(V)$ introduced by Y.-Z. Huang. 
Among the main results, we establish a bijection between the set of equivalence classes of 
irreducible $\N$-graded  $\phi$-coordinated $V$-modules and 
the set of isomorphism classes of irreducible $\widetilde{A}(V)$-modules. 
We also show that for a vertex operator algebra, rationality, regularity, and fusion rules
are independent of the choice of the conformal vector.
\end{abstract}

\section{Introduction}
The purpose of  this paper is two-fold, one of which
 is to establish an analogue of Zhu's $A(V)$-theory for what were called $\phi$-coordinated modules 
for weak quantum vertex algebras and the other is to study the dependence of certain fundamental notions on the choice of 
the conformal vector of a vertex operator algebra. 

In vertex operator algebra theory, an effective tool to determine irreducible modules 
for a given vertex operator algebra (VOA for short) $V$,
is the associative algebra $A(V)$ introduced by Zhu  (see \cite{zhu1,zhu2}), 
which is now commonly referred to as the Zhu algebra.
The core of Zhu's $A(V)$-theory is the natural bijection between the set of equivalence classes of
irreducible $\N$-graded $V$-modules and the set of  isomorphism classes of irreducible $A(V)$-modules. 
This theory was extended further by Frenkel and Zhu (see \cite{fz}; cf. \cite{li-fusion})  
to obtain a theorem for determining fusion rules 
in terms of $A(V)$-bimodules.
As the first application, Zhu's algebra was used to determine irreducible modules
for the simple affine VOAs of positive integer levels by Frenkel-Zhu (see \cite{fz}) and 
 for the simple Virasoro VOAs 
of central charges in the minimal series by Dong-Mason-Zhu and by Wang (see \cite{dmz}, \cite{wang}).
Since then, Zhu algebra $A(V)$ and Frenkel-Zhu $A(V)$-bimodule have been generalized in various directions 
from vertex operator algebras to vertex operator superalgebras, 
from modules to twisted modules, and then to logarithmic twisted modules  
(cf. \cite{kw}, \cite{dlm-twisted, dlm-anv, dlm-twisted-algebra}, \cite{dj1,dj2}, \cite{hy}).
In particular, a sequence of associative algebras $A_{n}(V)$ for $n\in \N$ was associated  in \cite{dlm-anv} to
each VOA $V$ and it was proved that $V$ is rational if and only if
$A_{n}(V)$ for $n\in \N$ are finite-dimensional and semi-simple. 

In \cite{huang}, among other results  Huang introduced another associative algebra $\widetilde{A}(V)$ 
for each vertex operator algebra $V$ 
and  established a bijection between the set of equivalence classes of irreducible $\N$-graded $V$-modules 
and the set of  isomorphism classes of irreducible $\widetilde{A}(V)$-modules. He also proved that 
$\widetilde{A}(V)$ is in fact isomorphic to the Zhu algebra $A(V)$. 

In a series of papers  (see \cite{li-nonlocal, li-const}), we  developed 
a theory of (weak) quantum vertex algebras, where
the notion of (weak) quantum vertex algebra generalizes the notions of vertex super-algebra
and vertex color-algebra (see \cite{xu1}, \cite{xu2}) in a certain natural way. 
The introduction of this theory was inspired by Etingof-Kazhdan's theory 
of quantum vertex operator algebras (see \cite{EK}) and 
the main goal  is to establish a theory of quantum vertex algebras which can be canonically associated to
quantum affine algebras. 
Note that the structure of a quantum vertex algebra is governed by a rational quantum Yang-Baxter operator 
whereas the structure of a quantum affine algebra is governed by a trigonometric quantum Yang-Baxter operator.
In order to associate quantum vertex algebras to quantum affine algebras, 
we furthermore developed a theory of what we called $\phi$-coordinated quasi modules 
for weak quantum vertex algebras  (see \cite{li-phi, li-equiv}). 
In this theory, $\phi$ stands for what we called an associate of the one-dimensional additive formal group (law) 
$F_{\add}(x,y)=x+y$. 
Indeed, by using this very theory with $\phi(x,z)=xe^{z}$ (a particular associate), weak quantum vertex algebras
were successfully  associated to quantum affine algebras {\em conceptually}.

As it was pointed out by Borcherds  (see \cite{Bo}), the formal group law $F_{\add}(x,y)$ underlies
the theory of vertex  algebras. This is apparent from the associativity
$$Y(u,z+x)Y(v,x)=Y(Y(u,z)v,x)$$
(an unrigorous version). 
This is also the case for the theory of (weak) quantum vertex algebras and modules, in which the same associativity 
is postulated though the usual locality (commutativity) is generalized to a braided locality. 

An {\em associate} of $F_{\add}(x,y)$ by definition (see \cite{li-phi}) is  a formal series $\phi(x,z)\in \C((x))[[z]]$, 
satisfying the conditions
\begin{eqnarray*}
\phi(x,0)=x,\    \   \    \   \phi(\phi(x,y),z)=\phi(x,y+z)\ \left(=\phi(x,F_{\add}(y,z))\right).
\end{eqnarray*}
It was proved that for each $p(x)\in \C((x))$, $e^{zp(x)\frac{d}{dx}}(x)$ is an associate 
of $F_{\add}(x,y)$ and every associate is of this form.
In particular, taking $p(x)=1$ and $x$, we obtain particular associates
\begin{eqnarray}
\phi(x,z)=e^{z\frac{d}{dx}}(x)=x+z=F_{\add}(x,z)\   \mbox{ and }\  \phi(x,z)=e^{zx\frac{d}{dx}}(x)= xe^{z}.
\end{eqnarray}
The essence  is that  to each associate $\phi(x,z)$, a theory of $\phi$-coordinated modules
for any weak quantum vertex algebra $V$ was attached, where an unrigorous version of
the defining associativity of a $\phi$-coordinated $V$-module $(W,Y_{W})$ is
 \begin{eqnarray}
Y_{W}(u,x_1)Y_{W}(v,x)|_{x_1=\phi(x,z)}=Y_{W}(Y(u,z)v,x).
\end{eqnarray}

For vertex operator algebras, there is a canonical connection between the notion of $\phi$-coordinated module and that of module.
Let $V$ be a VOA and $(W,Y_{W})$ a $V$-module.
 For $v\in V$, set $X_{W}(v,x)=Y_{W}(x^{L(0)}v,x)$.
 It was shown in \cite{li-vfa} that $(W,X_{W})$ carries the structure of a
 $\phi$-coordinated module with $\phi(x,z)=xe^{z}$ for the 
new VOA  obtained by Zhu (see \cite{zhu1, zhu2}, cf. \cite{lep0}) on $V$, where the new vertex operator map 
denoted by $Y[v,x]$ is given by $Y[v,x]=Y(e^{xL(0)},e^{x}-1)$ for $v\in V$. 
(In fact, one can show that this gives rise to an isomorphism of categories.)
This can be considered as a reformulation of a previous result of Lepowsky (see \cite{lep1}). 
It is also closely related to some of Huang's results (see \cite{huang}). 

Note that a VOA by definition is a vertex algebra $V$ equipped with
a conformal vector $\omega$, which provides a representation of the Virasoro algebra and 
a canonical $\Z$-grading on $V$. For a VOA $V$, 
the fundamental notions of module, rationality, regularity, and intertwining operator depend
on the conformal vector. A natural question is how much  these notions depend on the choice of the conformal vector.
 
In this paper, we exhibit a canonical connection between $\phi$-coordinated modules for a weak quantum vertex algebra
$V$ and the associative algebra $\widetilde{A}(V)$ introduced by Huang. 
Among the main results, to a weak quantum vertex algebra $V$ we associate a sequence of associative algebras 
$\widetilde{A}_{n}(V)$ for $n\in \N$ with $\widetilde{A}_{0}(V)=\widetilde{A}(V)$ and we establish a bijection between the set of equivalence classes of irreducible $\N$-graded  $\phi$-coordinated $V$-modules and 
the set of isomorphism classes of irreducible $\widetilde{A}(V)$-modules. 
On the other hand, we show that if a vertex algebra $V$ is a vertex operator algebra 
for some  conformal vector,
then rationality, regularity, and fusion rules are independent of the choice of the conformal vector. 

In the following, we describe the main content of this paper with some technical details.
First of all, recall that weak quantum vertex algebras are generalizations of vertex super-algebras.
Let $V$ be a weak quantum vertex algebra. 
The major axiom in the definition of a $\phi$-coordinated $V$-module $W$
with the vertex operator map $Y_{W}(\cdot,x)$ is that for $u,v\in V$,
there exists $k\in \N$ such that
$$(x_1-x_2)^{k}Y_{W}(u,x_1)Y_{W}(v,x_2)\in \Hom (W,W((x_1,x_2))),$$
$$x_2^{k}(e^{z}-1)^{k}Y_{W}(Y(u,z)v,x_2)=\left((x_1-x_2)^{k}Y_{W}(u,x_1)Y_{W}(v,x_2)\right)|_{x_1=x_2e^z}.$$
A $\Z$-graded (resp. $\N$-graded) $\phi$-coordinated $V$-module  
is a $\phi$-coordinated $V$-module $W$ equipped with a $\Z$-grading $W=\oplus_{n\in \Z}W[n]$ such that
\begin{eqnarray}
 v_{[m]} W[n]\subset W[m+n]\   \   \   \mbox{ for }v\in V,\ m,n\in \Z,
\end{eqnarray}
 where $Y_{W}(v,x)=\sum_{n\in \Z}v_{[n]}x^{-n}$.
Notice that the notion of $\Z$-graded $\phi$-coordinated module does {\em not} involve any grading on the ``algebra,''
unlike the notion of $\Z$-graded module which is defined for a $\Z$-graded vertex algebra. 

Let $V$ be a weak quantum vertex algebra
with a constant $\S$-map (see Section 2 for the definition),  e.g., a vertex super-algebra.
In this paper,  to $V$ we associate a sequence of associative algebras 
$\widetilde{A}_{n}(V)$ for $n\in \N$ with $\widetilde{A}(V)=\widetilde{A}_{0}(V)$, which 
are analogues of associative algebras $A_{n}(V)$ introduced in \cite{dlm-anv} 
for a vertex operator algebra $V$. As a vector space, $\widetilde{A}_{n}(V)=V/\widetilde{O}_{n}(V)$, where 
$\widetilde{O}_{n}(V)$ is the subspace linearly spanned by $\D V$ and by vectors
$$u\tilde{\circ}_{n}v=\Res_{x}\frac{e^{(n+1)x}}{(e^x-1)^{2n+2}}Y(u,x)v\   \   \mbox{ for }u,v\in V.$$
We show that if $W=\oplus_{m\in \N}W[m]$ is a (resp. irreducible) $\N$-graded 
$\phi$-coordinated $V$-module with $W[0]\ne 0$, then $W[0]$ is a (resp. irreducible) $\widetilde{A}(V)$-module
with $v$ acting as $v_{[0]}$ for $v\in V$.
 On the other hand,  we associate an associative algebra $\widetilde{U}(V)$ to $V$ and  
under a certain assumption on $\widetilde{U}(V)$, for any $\widetilde{A}(V)$-module $U$, we construct
an $\N$-graded $\phi$-coordinated $V$-module $W=\oplus_{n\in \N}W[n]$ with $W[0]\simeq U$ as an $\widetilde{A}(V)$-module.
In the case that $V$ is a vertex super-algebra,  we associate a Lie super-algebra to $V$ and then by using the P-B-W theorem 
we establish the validity of the assumption on $\widetilde{U}(V)$. 
Therefore,  for a general vertex super-algebra $V$ we obtain a bijection between the set of equivalence classes of
 irreducible $\N$-graded $\phi$-coordinated $V$-modules and the set of  isomorphism classes of irreducible 
 $\widetilde{A}(V)$-modules.  We believe that the assumption holds true for a more general case.
 We plan to study this problem in a separate publication.

Let $(V,\omega)$ denote a vertex operator algebra where $V$ is a vertex algebra and $\omega$ is a conformal vector.
 We denote by $\exp(V,\omega)$ Zhu's new vertex operator algebra associated to $(V,\omega)$. 
It was proved by Zhu (see \cite{zhu1,zhu2}, \cite{huang1, huang2}) that $\exp (V,\omega)$ is isomorphic to $V$.
Here, we show that $A_{n}(V,\omega)=\widetilde{A}_{n}(\exp(V,\omega))$ for $n\in \N$.
These enable us to show that 
$A_{n}(V,\omega)\simeq A_{n}(V,\omega')$ for $n\in \N$, where $\omega'$ is any conformal vector.
Then by applying the result of \cite{dlm-anv}  we show that the notions of rationality and regularity are 
independent of the choice of the conformal vector. 

Let $(V,\omega)$ be a VOA  and
let $W_1,W_2,W_3$ be modules for $V$ viewed as a vertex algebra.
From the definition of an intertwining operator ${\mathcal{Y}}(\cdot,x)$
of type $\binom{W_3}{W_1 W_2}$  (see \cite{fhl}), we see that the only axiom which uses the conformal vector is
$${\mathcal{Y}}(L(-1)w_1,x)=\frac{d}{dx}{\mathcal{Y}}(w_1,x)\   \   \mbox{ for }w_1\in W_1.$$
Assume that $W_1=\oplus_{n\in \N}W_{1}(n)$ is an $\N$-graded weak $(V,\omega)$-module such that
$W_1(0)$ is an irreducible $A(V)$-module and $W_1(0)$ generates $W_1$ as a $V$-module.
Let $\omega'$ be any conformal vector of $V$. 
We show that the spaces of intertwining operators of type $\binom{W_3}{W_1 W_2}$ with respect to $\omega$ 
and with respect to $\omega'$ are canonically isomorphic. 
 
This paper is organized as follows: In Section 2, we associate a sequence of associative algebras 
$\widetilde{A}_{n}(V)$ for $n\in \N$ to a weak quantum vertex algebra $V$ and associate an
$\widetilde{A}_{0}(V)$-module to each $\phi$-coordinated $V$-module.  In Section 3, from an 
$\widetilde{A}_{0}(V)$-module $U$ we construct an $\N$-graded $\phi$-coordinated $V$-module $W$ with 
$W[0]\simeq U$. In Section 4, we study the dependence of rationality, regularity and fusion rule on the 
choice of the conformal vector.

{\bf Acknowledgement.}
 We would like to thank Yufeng Pei for extensive discussions.
This research was partially supported by China NSF grants (Nos. 11471268, 11571391).

\section{$\phi$-coordinated modules for quantum vertex algebras and associative algebras $\widetilde{A}_{n}(V)$ }

In this section, first we recall the definitions of a weak quantum vertex algebra and a $\phi$-coordinated module for a 
weak quantum vertex algebra, and then we associate associative algebras $\widetilde{A}_{n}(V)$ for $n\in \N$ to each
weak quantum vertex algebra $V$ with a constant $\S$-map and we associate an $\widetilde{A}_{n}(V)$-module
$\widetilde{\Omega}_{n}(W)$ to every $\phi$-coordinated module $W$.

We begin by recalling from  \cite{li-nonlocal} the notion of weak quantum vertex algebra.
A {\em weak quantum vertex algebra}  is a vector space $V$ equipped with a linear map
\begin{eqnarray}
Y(\cdot,x):&& V\rightarrow \Hom (V,V((x)))\subset (\End V)[[x,x^{-1}]],\nonumber\\
&&v\mapsto Y(v,x)=\sum_{n\in \Z}v_{n}x^{-n-1}\   \   \left(\mbox{where }v_{n}\in \End V\right)
\end{eqnarray}
and with a vector ${\bf 1}$, called the {\em vacuum vector}, satisfying the conditions that 
$$Y({\bf 1},x)v=v,\   \   Y(v,x){\bf 1}\in V[[x]]\  \mbox{ and }\ \lim_{x\rightarrow 0}Y(v,x){\bf 1}=v\   \   \mbox{ for }v\in V$$
and that for $u,v\in V$, there exist (finitely many)
$$u^{(i)},v^{(i)}\in V,\  f_{i}(x)\in \C((x))\  \ (i=1,\dots,r)$$ 
such that
\begin{eqnarray}\label{eS-jacobi}
&&x_0^{-1}\delta\left(\frac{x_1-x_2}{x_0}\right)Y(u,x_1)Y(v,x_2)\nonumber\\
&&\hspace{1cm}-x_0^{-1}\delta\left(\frac{x_2-x_1}{-x_0}\right)\sum_{i=1}^{r}f_{i}(x_2-x_1)Y(v^{(i)},x_2)Y(u^{(i)},x_1)\nonumber\\
&&=x_2^{-1}\delta\left(\frac{x_1-x_0}{x_2}\right)Y(Y(u,x_0)v,x_2).
\end{eqnarray}

Let $V$ be a weak quantum vertex algebra. Define a linear operator $\D$ on $V$ by
$$\D (v)=v_{-2}{\bf 1}=\lim_{x\rightarrow 0}\frac{d}{dx}Y(v,x){\bf 1}
\   \   \mbox{ for }v\in V.$$ 
Then
\begin{eqnarray}
[\D,Y(v,x)]=Y(\D v, x)=\frac{d}{dx}Y(v,x).
\end{eqnarray}

The following notion was introduced in  \cite{li-phi}:

\bd{dphi-module}
{\em Let $V$ be a weak quantum vertex algebra. 
A {\em $\phi$-coordinated $V$-module} is a vector space $W$ 
equipped with a linear map
\begin{eqnarray}
Y_{W}(\cdot,x):&& V\rightarrow  (\End W)[[x,x^{-1}]],\nonumber\\
&&v\mapsto Y_{W}(v,x)=\sum_{n\in \Z}v_{[n]}x^{-n}\   \   \left(\mbox{where }v_{[n]}\in \End W\right),
\end{eqnarray}
satisfying the conditions that 
$$Y_{W}(v,x)w\in W((x))\    \    \    \mbox{ for }v\in V,\ w\in W,$$
$$Y_{W}({\bf 1},x)=1_{W},$$
and that for $u,v\in V$, there exists a nonnegative integer $k$ such that
$$(x_1-x_2)^{k}Y_{W}(u,x_1)Y_{W}(v,x_2)\in \Hom (W,W((x_1,x_2))),$$
\begin{eqnarray}\label{ephi-assoc}
(e^{x_0}-1)^{k}Y_{W}(Y(u,x_0)v,x_2)=\left((x_1/x_2-1)^{k}Y_{W}(u,x_1)Y_{W}(v,x_2)\right)|_{x_1=x_2e^{x_0}}.
\end{eqnarray}
We refer the last condition above as the {\em compatibility and weak $\phi$-associativity}. }
\ed

Let $\C(x)$ denote the field of rational functions.
For any $g(x)\in \C(x)$, we denote by $\tilde{g}(x)$ the formal Laurent series expansion of $g(x)$ at $0$.
On the other hand, note that  $p(e^x)\in \C[[x]]$ for any $p(x)\in \C[x]$ and $p(x)\mapsto p(e^x)$ is an algebra 
 embedding of $\C[x]$ into $\C[[x]]$, which gives rise to 
an algebra (field) embedding of $\C(x)$ into $\C((x))$. As a convention, for $g(x)\in \C(x)$,  we 
use $g(e^{x})$ to denote the image of $g(x)$ in $\C((x))$ under this embedding.

The following two results were obtained also in \cite{li-phi}:

\bl{lD-property}
Let $(W,Y_{W})$ be any $\phi$-coordinated $V$-module. Then
\begin{eqnarray}\label{ederivative-property}
Y_{W}(\D v,x)=\left(x\frac{d}{dx}\right)Y_{W}(v,x)
\end{eqnarray}
for $v\in V$. 
\el

\bp{pphi-jacobi}
Let $(W,Y_{W})$ be a $\phi$-coordinated $V$-module and let  $u,v\in V$.  Suppose that there are
$u^{(i)},v^{(i)}\in V$ and $f_{i}(x)\in \C((x))$  $(i=1,\dots, r)$ such that (\ref{eS-jacobi}) holds, where $f_i(x)=g_i(e^{-x})$
with $g_i(x)\in \C(x)$. Then
\begin{eqnarray}\label{phi-Jacobi-qva-0}
&&(zx_2)^{-1} \delta\left(\frac{x_1-x_2}{zx_2}\right)Y_{W}(u, x_1)Y_{W}(v, x_2)\nonumber\\
&&\hspace{1cm}-(zx_2)^{-1}\delta \left(\frac{x_2-x_1}{-zx_2}\right)\sum_{i=1}^{r}\tilde{g_{i}}(x_1/x_2)
Y_{W}(v^{(i)}, x_2)Y_{W}(u^{(i)}, x_1) \nonumber\\
&&=x_1^{-1}\delta\left(\frac{x_2(1+z)}{x_1}\right)Y_{W}(Y(u, \log(1 + z))v, x_2),
\end{eqnarray}
where
\begin{eqnarray*}
\log (1+z)=\sum_{n\ge 1}(-1)^{n-1}\frac{1}{n}z^{n}\in \C[[z]].
\end{eqnarray*}
\ep

\br{rsub-existence}
{\em Let $A(x), B(x)\in \Hom (W,W((x)))$ with $W$ a vector space. We mention that for any $w\in W$, 
$$A((z+1)x)B(x)w\   \mbox{ exists in }W[[z,z^{-1},x,x^{-1}]].$$
Indeed, writing $A(x)=\sum_{n\in \Z}a(n)x^{-n},\  B(x)=\sum_{n\in \Z}b(n)x^{-n}$, we have 
\begin{eqnarray*}
A((z+1)x)B(x)w=\sum_{m,n\in \Z}\sum_{j\ge 0}\binom{-m}{j}z^{-m-j}x^{-m-n}a(m)b(n)w.
\end{eqnarray*}
Let $s\in \Z$ be such that $b(n)w=0$ for $n>s$.
For any integers $p,q$, the coefficient of $z^{p}x_2^{q}$ in $A((z+1)x)B(x)w$
is the sum of terms $\binom{-m}{j}a(m)b(n)w$ with
$m,n\in \Z,\ j\in \N$, satisfying the relations
$$m+n=-q,\   \  n\le s,\  \   -m-j=p,$$
which have only finitely many solutions.}
\er

\bp{pweak-assoc}
In the definition of a $\phi$-coordinated $V$-module, the last axiom can be equivalently replaced by 
{\em weak $\phi$-associativity:} For any $u,v\in V,\ w\in W$, there exists $l\in \N$ such that
\begin{eqnarray}
(z+1)^{l}Y_{W}(u,(z+1)x_2)Y_{W}(v,x_2)w=(1+z)^{l}Y_{W}(Y(u,\log (1+z))v,x_2)w.
\end{eqnarray}
\ep

\begin{proof}  Suppose that $(W,Y_{W})$ is a $\phi$-coordinated $V$-module.  Let $u,v\in V,\ w\in W$.
By Proposition \ref{pphi-jacobi}, Jacobi identity relation (\ref{phi-Jacobi-qva}) holds. Let $l\in \N$ be such that
$$x_1^{l}\tilde{g_{i}}(x_1/x_2)Y_{W}(u^{(i)},x_1)w\in W[[x_1,x_2,x_2^{-1}]]\   \   \   \mbox{ for }i=1,\dots,r.$$
Then applying $\Res_{x_1}x_1^{l}$ to (\ref{phi-Jacobi-qva}) which is applied to vector $w$,
we obtain
\begin{eqnarray*}
(z+1)^{l}x_2 ^{l}Y_{W}(u, (z+1)x_2)Y_{W}(v,x_2)w=(z+1)^{l}x_2^{l}Y_{W}(Y(u,\log(1+z))v,x_2)w.
\end{eqnarray*}
This proves that weak $\phi$-associativity holds. 
 
On the other hand, assume weak $\phi$-associativity. Let $u,v\in V$.
There exists $k\in \N$ such that $x^{k}Y(u,x)v\in V$.  For any $w\in W$, there exists $l\in \N$ such that
\begin{eqnarray*}
&&(z+1)^{l}(zx_2)^{k}Y_{W}(u, (z+1)x_2)Y_{W}(v,x_2)w\\
&=&(z+1)^{l}(zx_2)^{k}Y_{W}(Y(u,\log(1+z))v,x_2)w.
\end{eqnarray*}
We see that the right-hand side involves only nonnegative powers of $z$,
 and so does the left-hand side. This implies that
  $$x_1^{l}(x_1-x_2)^{k}Y_{W}(u,x_1)Y_{W}(v,x_2)w\in W[[x_1,x_2]][x_2^{-1}].$$
  As $k$ is independent of $w$, we have
 $$(x_1-x_2)^{k}Y_{W}(u,x_1)Y_{W}(v,x_2)\in \Hom(W,W((x_1,x_2))).$$
Considering the substitution $x_1=x_2e^{x_0}$ as the composition of substitutions $x_1=zx_2+x_2$ with $z=e^{x_0}-1$,
we obtain 
 \begin{eqnarray*}
&&(x_2e^{x_0})^{l}\left((x_1-x_2) ^{k}Y_{W}(u, x_1)Y_{W}(v,x_2)w\right)|_{x_1=x_2e^{x_0}}\\
&=&\left(\left(x_1^{l}(x_1-x_2) ^{k}Y_{W}(u, x_1)Y_{W}(v,x_2)w\right)|_{x_1=zx_2+x_2}\right)|_{z=e^{x_0}-1}\\
&=&\left((zx_2+x_2)^{l}(zx_2)^{k}Y_{W}(Y(u,\log(1+z))v,x_2)w\right)|_{z=e^{x_0}-1}\\
&=&(x_2e^{x_0})^{l}(x_2e^{x_0}-x_2)^{k}Y_{W}(Y(u,x_0)v,x_2)w.
\end{eqnarray*}
This proves that weak $\phi$-associativity holds.
\end{proof}

Let $(W,Y_{W})$ be a $\phi$-coordinated $V$-module and let  $u,v\in V$.  By applying $\Res_{z}x_2$ to both sides 
of (\ref{phi-Jacobi-qva}) we get
\begin{eqnarray}\label{phi-Scommutator}
&&Y_{W}(u, x_1)Y_{W}(v, x_2)-\sum_{i=1}^{r}\tilde{g_{i}}(x_1/x_2)
Y_{W}(v^{(i)}, x_2)Y_{W}(u^{(i)}, x_1) \nonumber\\
&=&\Res_{z}x_2
x_1^{-1}\delta\left(\frac{x_2(1+z)}{x_1}\right)Y_{W}(Y(u, \log(1 + z))v, x_2)\nonumber\\
&=&\Res_{x_0}\delta\left(\frac{x_2e^{x_0}}{x_1}\right)Y_{W}(Y(u,x_0)v,x_2).
\end{eqnarray}
In terms of components we have
\begin{eqnarray}\label{need}
u_{[m]}v_{[n]}-\sum_{i=1}^{r}\sum_{k\in \Z}a_{i,k}v^{(i)}_{[n-k]}u^{(i)}_{[m+k]}=\sum_{j\ge 0}\frac{m^{j}}{j!} (u_{j}v)_{[m+n]},
\end{eqnarray}
where $\tilde{g_{i}}(x)=\sum_{k\in \Z}a_{i,k}x^{k}$. Notice that the relation above is homogeneous 
with $\deg a_{[n]}=n$ for $a\in V,\ n\in \Z$. 

\bd{dgraded-modules}
{\em A {\em $\Z$-graded (resp. $\N$-graded) $\phi$-coordinated $V$-module}
 is a $\phi$-coordinated $V$-module $W$ equipped with a $\Z$-grading $W=\oplus_{n\in \Z}W_{[n]}$ 
(resp. $\N$-grading) such that 
\begin{eqnarray}
v_{[m]}W_{[n]}\subset W_{[m+n]}\   \   \   \mbox{ for }v\in V,\ m,n\in \Z.
\end{eqnarray}}
\ed

To associate associative algebras to a weak quantum vertex algebra $V$, we shall need to assume
that $V$ is of a certain type.

\bd{dwqva-constant}
{\em A weak quantum vertex algebra $V$ is called a {\em weak quantum vertex algebra with a constant $\S$-map} if
for any $u,v\in V$, there exist $u^{(i)},v^{(i)}\in V,\ 1\le i\le r$ such that
\begin{eqnarray}\label{e2.12}
&&x_0^{-1}\delta\left(\frac{x_1-x_2}{x_0}\right)Y(u,x_1)Y(v,x_2)
-x_0^{-1}\delta\left(\frac{x_2-x_1}{-x_0}\right)\sum_{i=1}^{r}Y(v^{(i)},x_2)Y(u^{(i)},x_1)\nonumber\\
&&\hspace{3cm}=x_2^{-1}\delta\left(\frac{x_1-x_0}{x_2}\right)Y(Y(u,x_0)v,x_2).
\end{eqnarray}}
\ed

Note that this family of weak quantum vertex algebras includes vertex algebras, vertex super-algebras, 
and vertex color-algebras (see \cite{xu1, xu2}).
Let $V$ be a weak quantum vertex algebra with a constant $\S$-map and let $(W,Y_{W})$ be a $\phi$-coordinated $V$-module.
For $u,v\in V$, there exist vectors $u^{(i)},v^{(i)}\in V$ $(i=1,\dots,r)$ such that
\begin{eqnarray}\label{phi-Jacobi-qva}
&&(zx_2)^{-1} \delta\left(\frac{x_1-x_2}{zx_2}\right)Y_{W}(u, x_1)Y_{W}(v, x_2)\nonumber\\
&&\hspace{1cm}-(zx_2)^{-1}\delta \left(\frac{x_2-x_1}{-zx_2}\right)\sum_{i=1}^{r}
Y_{W}(v^{(i)}, x_2)Y_{W}(u^{(i)}, x_1) \nonumber\\
&&=x_1^{-1}\delta\left(\frac{x_2(1+z)}{x_1}\right)Y_{W}(Y(u, \log(1 + z))v, x_2).
\end{eqnarray}
In particular, we have
\begin{eqnarray}\label{need-2}
u_{[m]}v_{[n]}-\sum_{i=1}^{r}v^{(i)}_{[n]}u^{(i)}_{[m]}=\sum_{j\ge 0}\frac{m^{j}}{j!} (u_{j}v)_{[m+n]}.
\end{eqnarray}

We recall the associative algebra $\widetilde{A}(V)$ introduced by Huang. 
Let $V$ be a vertex algebra. Follow \cite{huang} to define a bilinear operation $*$ on $V$ by
\begin{eqnarray}
u*v=\Res_{x}\frac{e^{x}}{e^{x}-1}Y(u,x)v=\Res_{z}\frac{1}{z}Y(u,\log (1+z))v
\end{eqnarray}
for $u,v\in V$, where $\frac{1}{e^{x}-1}$ is considered as the inverse of $e^{x}-1$ in $\C((x))$.  
Let $\widetilde{O}(V)$ be the linear span of
\begin{eqnarray}
\Res_{x}\frac{e^{x}}{(e^{x}-1)^{2}}Y(u,x)v=\Res_{z}\frac{1}{z^{2}}Y(u,\log (1+z))v
\end{eqnarray}
for $u,v\in V$. The following was proved in \cite{huang}:

\bt{thuang}
Let $V$ be a vertex algebra. Then $\widetilde{O}(V)$ is a two-sided ideal of $(V,*)$ and
the quotient $V/\widetilde{O}(V)$, denoted by $\widetilde{A}(V)$, is an associative algebra with ${\bf 1}+\widetilde{O}(V)$ 
as an identity.
\et

\br{rhuang}
{\em Assuming that $V$ is a vertex operator algebra in the sense of \cite{flm} and \cite{fhl},
Huang established in \cite{huang} a one-to-one correspondence
between the set of isomorphism classes of irreducible $\widetilde{A}(V)$-modules and the set of equivalence classes of 
irreducible $\N$-graded $V$-modules. Huang also proved that 
$\widetilde{A}(V)$ is in fact isomorphic to the Zhu algebra $A(V)$. }
\er

Next, we generalize Huang's result from a vertex algebra to a weak quantum vertex algebra with a constant $\S$-map
in a more general setting.

\bd{dA(V)}
{\em Let $V$ be a weak quantum vertex algebra. For any $f(x)\in \C((x))$, we define a bilinear operation $*_{f(x)}$ on $V$ by
\begin{eqnarray}
u*_{f(x)}v=\Res_{x}f(x)Y(u,x)v
\end{eqnarray}
for $u,v\in V$.  Furthermore, for any subspace $B$ of $\C((x))$ set
\begin{eqnarray}
O_{B}(V)={\rm span}\{ u*_{g(x)}v\ |\  u,v\in V,\ g(x)\in B\}.
\end{eqnarray}}
\ed

\bp{pap-algebra}
Let $V$ be a weak quantum vertex algebra with a constant  $\S$-map.
Let $f(x)\in \C((x))$  and let $B$ be a subspace of $\C((x))$ such that
\begin{eqnarray}\label{eBconditions}
f^{(j)}(x)g(x),\    \    f(x)g^{(j)}(-x),\   \   f(x)f^{(j)}(-x),\   \  f^{(j+1)}(x)\in B
\end{eqnarray}
for $g(x)\in B$, $j\ge 0$, where $f^{(j)}(x)$ denotes the formal $j$-th derivative.
Then the subspace $O_{B}(V)$ of $V$ is a two-sided ideal of
the non-associative algebra $(V,*_{f(x)})$. Furthermore, $V/O_{B}(V)$ is an associative algebra. If 
$f(x)\in x^{-1}\C[[x]]$ with $\Res_{x}f(x)=1$,
then ${\bf 1}+O_{B}(V)$ is an identity element.
\ep

\begin{proof} It essentially follows from Zhu's arguments in \cite{zhu1}.
First of all, for $g(x)\in B,\ u,v\in V,\ j\in \N$, we have
$$\Res_{x}g(x)Y(\D^{j}u,x)v=\Res_{x}g(x)\left(\frac{d}{dx}\right)^{j}Y(u,x)v
=(-1)^{j}\Res_{x}g^{(j)}(x)Y(u,x)v.$$
Then
\begin{eqnarray}\label{efirst-fact}
\Res_{x}g^{(j)}(x)Y(u,x)v\in O_{B}(V)
\end{eqnarray}
for all $g(x)\in B,\ j\in \N,\ u,v\in V$.

For convenience, in the following we shall simply write $u*v$ for $u*_{f(x)}v$.
Let $u,v,w\in V,\ g(x)\in B$. We have
\begin{eqnarray*}
&&u*(v*_{g} w)\\
&=&\Res_{x_{1}}\Res_{x_{2}}f(x_{1})g(x_{2})Y(u,x_{1})Y(v,x_{2})w\\
&=&\Res_{x_{1}}\Res_{x_{2}}f(x_{1})g(x_{2})\sum_{i=1}^{r}Y(v^{(i)},x_{2})Y(u^{(i)},x_{1})w\\
&&-\sum_{j\ge 0}\Res_{x_{1}}\Res_{x_{2}}f(x_{1})g(x_{2})
Y(u_{j}v,x_{2})\frac{1}{j!}\left(\frac{\partial}{\partial x_{2}}\right)^{j}x_{1}^{-1}\delta\left(\frac{x_{2}}{x_{1}}\right)\\
&=&\sum_{i=1}^{r}\Res_{x_{2}}g(x_{2})Y(v^{(i)},x_{2})\left(\Res_{x_{1}}f(x_{1})Y(u^{(i)},x_{1})w\right)\\
&&-\sum_{j\ge 0}\Res_{x_{1}}\Res_{x_{2}}g(x_{2})
Y(u_{j}v,x_{2})w\frac{1}{j!}\left(\frac{\partial}{\partial x_{2}}\right)^{j}\left[x_{1}^{-1}\delta\left(\frac{x_{2}}{x_{1}}\right)f(x_{2})\right]\\
&=&\sum_{i=1}^{r}v^{(i)}*_{g} (u^{(i)}*w)-\sum_{j\ge 0}\frac{1}{j!}\Res_{x_{2}}g(x_{2})f^{(j)}(x_{2})
Y(u_{j}v,x_{2})w,
\end{eqnarray*}
\begin{eqnarray*}
&&(u*_{g} v)*w\\
&=&\Res_{x_{2}}\Res_{x_{0}}f(x_{2})g(x_{0})Y(Y(u,x_{0})v,x_{2})w\\
&=&\Res_{x_{2}}\Res_{x_{1}}f(x_{2})g(x_{1}-x_{2})Y(u,x_{1})Y(v,x_{2})w\\
&&-\sum_{i=1}^{r}\Res_{x_{2}}\Res_{x_{1}}f(x_{2})g(-x_{2}+x_{1})Y(v^{(i)},x_{2})Y(u^{(i)},x_{1})w\\
&=&\Res_{x_{2}}\Res_{x_{1}}f(x_{2})\left(e^{-x_{2}\frac{\partial}{\partial x_{1}}}g(x_{1})\right)Y(u,x_{1})Y(v,x_{2})w\\
&&-\sum_{i=1}^{r}\Res_{x_{2}}\Res_{x_{1}}f(x_{2})\left(e^{-x_{1}\frac{\partial}{\partial x_{2}}}g(-x_{2})\right)Y(v^{(i)},x_{2})Y(u^{(i)},x_{1})w\\
&=&\sum_{j\ge 0}\frac{1}{j!}\Res_{x_{1}}g^{(j)}(x_{1})Y(u,x_{1})
\left(\Res_{x_{2}}(-x_{2})^{j}f(x_{2})Y(v,x_{2})w\right)\\
&&-\sum_{i=1}^{r}\sum_{j\ge 0}\frac{1}{j!}\Res_{x_{2}}f(x_{2}) g^{(j)}(-x_{2})Y(v^{(i)},x_{2})\left(\Res_{x_{1}}x_{1}^{j}Y(u^{(i)},x_{1})w\right).
\end{eqnarray*}
Then it follows from (\ref{eBconditions}) and (\ref{efirst-fact}) that $O_{B}(V)$ is a two-sided ideal of
the non-associative algebra $(V,*_{f(x)})$. Furthermore, for $u,v,w\in V$ we have
\begin{eqnarray*}
&&(u*v)*w\\
&=&\Res_{x_{2}}\Res_{x_{0}}f(x_{2})f(x_{0})Y(Y(u,x_{0})v,x_{2})w\\
&=&\Res_{x_{2}}\Res_{x_{1}}f(x_{2})f(x_{1}-x_{2})Y(u,x_{1})Y(v,x_{2})w\\
&&-\sum_{i=1}^{r}\Res_{x_{2}}\Res_{x_{1}}f(x_{2})f(-x_{2}+x_{1})Y(v^{(i)},x_{2})Y(u^{(i)},x_{1})w\\
&=&\Res_{x_{2}}\Res_{x_{1}}f(x_{2})\left(e^{-x_{2}\frac{\partial}{\partial x_{1}}}f(x_{1})\right)Y(u,x_{1})Y(v,x_{2})w\\
&&-\sum_{i=1}^{r}\Res_{x_{2}}\Res_{x_{1}}f(x_{2})\left(e^{-x_{1}\frac{\partial}{\partial x_{2}}}f(-x_{2})\right)Y(v^{(i)},x_{2})Y(u^{(i)},x_{1})w\\
&=&\Res_{x_{2}}\Res_{x_{1}}f(x_{2})f(x_{1})Y(u,x_{1})Y(v,x_{2})w\\
&&+\sum_{n\ge 1}\frac{1}{n!}\Res_{x_{1}}f^{(n)}(x_{1})Y(u,x_{1})
\left(\Res_{x_{2}}(-x_{2})^{n}f(x_{2})Y(v,x_{2})w\right)\\
&&-\sum_{i=1}^{r}\sum_{j\ge 0}\frac{1}{j!}\Res_{x_{2}}f(x_{2}) f^{(j)}(-x_{2})Y(v^{(i)},x_{2})\left(\Res_{x_{1}}x_{1}^{j}Y(u^{(i)},x_{1})w\right)\\
&=&u*(v*w)\\
&&+\sum_{n\ge 1}\frac{1}{n!}\Res_{x_{1}}f^{(n)}(x_{1})Y(u,x_{1})
\left(\Res_{x_{2}}(-x_{2})^{n}f(x_{2})Y(v,x_{2})w\right)\\
&&-\sum_{i=1}^{r}\sum_{j\ge 0}\frac{1}{j!}\Res_{x_{2}}f(x_{2}) f^{(j)}(-x_{2})Y(v^{(i)},x_{2})\left(\Res_{x_{1}}x_{1}^{j}Y(u^{(i)},x_{1})w\right).
\end{eqnarray*}
It follows from (\ref{eBconditions})  that $V/O_{B}(V)$ is an associative algebra. 

Now, assume that $f(x)\in x^{-1}\C[[x]]$ with $\Res_{x}f(x)=1$. For $v\in V$, we have
\begin{eqnarray*}
&&{\bf 1}*v=\Res_{x}f(x)Y({\bf 1},x)v=\Res_{x}f(x)v=v,\\
&& v*{\bf 1}=\Res_{x}f(x)Y(v,x){\bf 1}=\Res_{x}f(x)e^{x\D}v=v.
\end{eqnarray*}
Thus ${\bf 1}+O_{B}(V)$ is an identity of $V/O_{B}(V)$.
\end{proof}

The following is straightforward:

\bl{lspecial}
Let $f(x)\in \C((x))$. Set $B=\span \{ f^{(k)}(x)\  |\  k\ge 1\}$. 
Then the conditions in (\ref{eBconditions}) amount to that $B$ is a subalgebra of $\C((x))$ such that 
\begin{eqnarray}\label{efBconditions-2}
f(x)B\subset B,\    \   f(x)f(-x),\  f(x)f^{(k)}(-x)\in B\   \    \mbox{ for }k\ge 1.
\end{eqnarray}
Furthermore, if $f^{(1)}(-x)\in B$, then the condition (\ref{efBconditions-2}) amounts to $f(x)B\subset B$ and $f(x)f(-x)\in B$.
\el
 
We also have the following result:

\bl{lD-ideal}
Let $V$ be a weak quantum vertex algebra with a constant  $\S$-map and let
$f(x)\in \C((x)),\  B\subset \C((x))$ as in Proposition \ref{pap-algebra} such that  (\ref{eBconditions}) holds.
Then the image of $\D V$ in $V/O_{B}(V)$ is a two-sided ideal.
\el

\begin{proof}  For $u,v\in V$, we have
\begin{eqnarray*}
(\D u)*_{f(x)}v=\Res_{x}f(x)Y(\D u,x)v=\Res_{x}f(x)\frac{d}{dx}Y(u,x)v
=-\Res_{x}f'(x)Y(u,x)v
\end{eqnarray*}
and
\begin{eqnarray*}
u*_{f(x)}(\D v)&=&\Res_{x}f(x)Y(u,x)\D v\\
&=&\Res_{x}\left(f(x)\D Y(u,x)v-f(x)Y(\D u,x)v\right)\\
&=&\D (u*_{f(x)}v)-(\D u)*_{f(x)}v.
\end{eqnarray*}
It follows immediately that the image of $\D V$ in $V/O_{B}(V)$ is a  two-sided ideal.
\end{proof}

The following consideration is motivated by the construction of the algebra $A_{n}(V)$ (see  \cite{dlm-anv}, page 69). 
Let $n$ be a nonnegative integer.
Set
\begin{eqnarray}
f_{n}(x)=\sum_{i=0}^{n}\binom{-n-1}{i}\frac{e^{(n+1)x}}{(e^{x}-1)^{i+n+1}}\in \C((x)),
\end{eqnarray}
\begin{eqnarray}
g_{n}(x)=\frac{e^{(n+1)x}}{(e^{x}-1)^{2n+2}}\in \C((x)),
\end{eqnarray}
where for any nonnegative integer $k$, $1/(e^x-1)^{k}$  denotes the inverse of $(e^x-1)^{k}$ in $\C((x))$.
We also set
\begin{eqnarray}
B_{n}=\span\left\{ \frac{e^{(n+1)x}}{(e^{x}-1)^{2n+2+k}}\ | \ k\in \N\right\}.
\end{eqnarray}

\bl{lanv}
Let $n$ be a nonnegative integer. Then all the conditions in (\ref{eBconditions}) are satisfied 
 with $B=B_{n}$ and $f(x)=f_{n}(x)$. 
\el

\begin{proof}  We shall apply Lemma \ref{lspecial}. First of all, we have
\begin{eqnarray}\label{efn-derivative}
&&f_{n}'(x)\nonumber\\
&=&\sum_{i=0}^{n}\binom{-n-1}{i}\left((n+1)\frac{e^{(n+1)x}}{(e^{x}-1)^{i+n+1}}
-(i+n+1)\frac{e^{(n+1)x}}{(e^{x}-1)^{i+n+2}}e^{x}\right)\nonumber\\
&=&\sum_{i=0}^{n}\binom{-n-1}{i}\left((-i)\frac{e^{(n+1)x}}{(e^{x}-1)^{i+n+1}}
-(i+n+1)\frac{e^{(n+1)x}}{(e^{x}-1)^{i+n+2}}\right)\nonumber\\
&=&\sum_{i=1}^{n}\binom{-n-1}{i}(-i)\frac{e^{(n+1)x}}{(e^{x}-1)^{i+n+1}}
-\sum_{i=0}^{n-1}\binom{-n-1}{i}(i+n+1)\frac{e^{(n+1)x}}{(e^{x}-1)^{i+n+2}}\nonumber\\
&&\  \   +\binom{-n-1}{n}(2n+1)\frac{e^{(n+1)x}}{(e^{x}-1)^{2n+2}}\nonumber\\
&=&\sum_{i=1}^{n}\left(\binom{-n-1}{i}(-i)
-\binom{-n-1}{i-1}(i+n)\right)\frac{e^{(n+1)x}}{(e^{x}-1)^{i+n+1}} \nonumber\\
&&\  \   +\binom{-n-1}{n}(2n+1)\frac{e^{(n+1)x}}{(e^{x}-1)^{2n+2}}\nonumber\\
&=&\binom{-n-1}{n}(2n+1)\frac{e^{(n+1)x}}{(e^{x}-1)^{2n+2}}.
\end{eqnarray}
That is, $f_{n}'(x)=\binom{-n-1}{n}(2n+1)g_{n}(x)\in B_{n}$. 
Since 
\begin{eqnarray}
&&\frac{d}{dx}\frac{e^{(n+1)x}}{(e^{x}-1)^{2n+2+k}}\nonumber\\
&=&-(n+1+k)\frac{e^{(n+1)x}}{(e^{x}-1)^{2n+2+k}}-(2n+2+k)\frac{e^{(n+1)x}}{(e^{x}-1)^{2n+3+k}}
\end{eqnarray}
for $k\in \N$, it follows that $B_{n}$ is closed under $\frac{d}{dx}$ and
\begin{eqnarray}\label{eBn-new}
B_{n}=\span\left\{ g_{n}^{(i)}(x)\ |\ i\in \N\right\}=\span\{ f_{n}^{(i+1)}(x)\ |\  i\in \N\}.
\end{eqnarray}
 
For any nonnegative integers $m$ and $k$,  noticing that 
\begin{eqnarray*}
e^{mx}=(e^{x})^{m}=\sum_{i=0}^{m}\binom{m}{i}(e^{x}-1)^{m-i},
\end{eqnarray*}
we have
\begin{eqnarray}\label{this}
\frac{e^{mx}}{(e^{x}-1)^{m+k}}\in {\rm span}\left\{ \frac{1}{(e^{x}-1)^{k+j}}\ |\ j=0,\dots,m  \right\}.
\end{eqnarray}
It follows that $B_{n}$ is a subalgebra of $\C((x))$. 
We also have  $f_{n}(x)B_{n}\subset  B_{n}$
as
\begin{eqnarray*}
f_{n}(x)\cdot \frac{e^{(n+1)x}}{(e^{x}-1)^{2n+2+k}}
=\sum_{i=0}^{n}\binom{-n-1}{i}\frac{e^{(n+1)x}\cdot e^{(n+1)x}}{(e^{x}-1)^{i+3n+3+k}}  \in B_{n}
\end{eqnarray*}
by (\ref{this})  for $k\in \N$.
Using (\ref{this}) again  we get
\begin{eqnarray*}
f_{n}(x)f_{n}(-x)&=&\sum_{i,j=0}^{n}\binom{-n-1}{i}\binom{-n-1}{j} \frac{e^{(n+1)x}}{(e^{x}-1)^{i+n+1}}
\cdot \frac{e^{-(n+1)x}}{(e^{-x}-1)^{j+n+1}}\\
&=&\sum_{i,j=0}^{n}\binom{-n-1}{i}\binom{-n-1}{j}(-1)^{j+n+1}\frac{e^{(n+1)x}\cdot e^{jx}}{(e^{x}-1)^{i+j+2n+2}}\in B_{n}.
\end{eqnarray*}
On the other hand,  we have $g(-x)\in B_{n}$ for every $g(x)\in B_{n}$ as
\begin{eqnarray*}
\frac{e^{-(n+1)x}}{(e^{-x}-1)^{2n+2+k}}=(-1)^{k}\frac{e^{(n+1)x}\cdot e^{kx}}{(e^{x}-1)^{2n+2+k}}\in B_{n}
\end{eqnarray*}
for $k\in \N$.  Then by Lemma \ref{lspecial} all the conditions in (\ref{eBconditions}) are satisfied.
\end{proof}

As an immediate consequence of Proposition \ref{pap-algebra} and Lemma \ref{lanv} we have:

\bc{canV}
Let $n$ be a nonnegative integer and let $f_{n}(x)\in \C((x)), \ B_{n}\subset \C((x))$ be defined as before.
Then for any weak quantum vertex algebra $V$ with a constant $\S$-map, $O_{B_n}(V)$ is an ideal of
the non-associative algebra $(V,*_{f_n(x)})$ and  the quotient 
 $V/O_{B_n}(V)$ is an associative algebra with ${\bf 1}+O_{B_n}(V)$ as an identity element.
\ec

Recall from Lemma \ref{lD-ideal} that the image of $\D V$ is a two-sided ideal of the associative algebra $V/O_{B_{n}}(V)$.
We now define an associative algebra which corresponds to $A_{n}(V)$ introduced in \cite{dlm-anv} 
for a vertex operator algebra $V$.

\bd{dtildeA}
{\em Let $V$ be a weak quantum vertex algebra with a constant $\S$-map and let $n$ be a nonnegative integer. 
We define $\widetilde{A}_{n}(V)$
to be the quotient algebra of $V/O_{B_{n}}(V)$ modulo the image of $\D V$ (which is an ideal).  }
\ed

For $u,v\in V$, set
\begin{eqnarray}
u\tilde{*}_{n}v&=&u*_{f_{n}(x)}v=\Res_{x}f_{n}(x)Y(u,x)v,\\
u\tilde{\circ}_{n}v&=&u*_{g_{n}(x)}v=\Res_{x}g_{n}(x)Y(u,x)v.
\end{eqnarray}
We have
\begin{eqnarray}
u\tilde{*}_{n}v&=&\Res_{x}\sum_{i=0}^{n}\binom{-n-1}{i}
\frac{e^{(n+1)x}}{(e^{x}-1)^{n+1+i}}Y(u,x)v\nonumber\\
&=&\Res_{z}\sum_{i=0}^{n}\binom{-n-1}{i}
\frac{(1+z)^{n}}{z^{n+1+i}}Y(u,\log (1+z))v
\end{eqnarray}
and
\begin{eqnarray}
u\tilde{\circ}_{n}v=\Res_{x}\frac{e^{(n+1)x}}{(e^{x}-1)^{2n+2}}Y(u,x)v
=\Res_{z}\frac{(1+z)^{n}}{z^{2n+2}}Y(u,\log (1+z))v.
\end{eqnarray}
Note that it follows from (\ref{eBn-new}) and (\ref{efirst-fact}) that
$O_{B_n}(V)$ is linearly spanned by $u\tilde{\circ}_{n}v$ for $u,v\in V$. 

Next,  we associate a (resp. irreducible) $\widetilde{A}(V)$-module to each (resp. 
irreducible) $\N$-graded $\phi$-coordinated $V$-module. 
 
\bp{pomega}
Let $(W,Y_{W})$ be any $\phi$-coordinated $V$-module and let $n\in \N$. Set
\begin{eqnarray}
\widetilde{\Omega}_{n}(W)=\{ w\in W\ |\  v_{[m]}w=0\  \ \mbox{ for }v\in V,\ m\ge n+1\}.
\end{eqnarray}
Then $\widetilde{\Omega}_{n}(W)$ is an $\widetilde{A}_{n}(V)$-module with $v$ acting as $v_{[0]}$ for $v\in V$.
Furthermore, if $W=\oplus_{m\in \N}W[m]$ is an irreducible $\N$-graded $\phi$-coordinated $V$-module with $W[0]\ne 0$, then 
$\widetilde{\Omega}_{0}(W)=W[0]$ and $W[0]$ is an irreducible $\widetilde{A}(V)$-module.
\ep

\begin{proof} We first prove that $\widetilde{\Omega}(W)$ is closed under the action of $v_{[0]}$ for every $v\in V$. 
Let $u,v\in V,\ w\in \widetilde{\Omega}_{n}(W)$. Using (\ref{need}) (or (\ref{need-2})) we have
\begin{eqnarray*}
u_{[p]}v_{[q]}w-\sum_{i=1}^{r}v^{(i)}_{[q]}u^{(i)}_{[p]}w=\sum_{j\ge 0}\frac{m^{j}}{j!} (u_{j}v)_{[p+q]}w
\end{eqnarray*}
for $p,q\in \Z$. Then it follows that 
$u_{[p]}v_{[q]}w=0$ for $p\ge n+1$ and $q\ge 0$. Thus $v_{[q]}w\in \widetilde{\Omega}_{n}(W)$ for $q\ge 0$. 
In particular, we have $v_{[0]}w\in \widetilde{\Omega}_{n}(W)$. 

Let $u,v\in V,\ w\in \widetilde{\Omega}_{n}(W)$.
With $w\in \widetilde{\Omega}_{n}(W)$, we have $x_1^{n}Y_{W}(u^{(i)},x_1)w\in W[[x_1]]$ for $1\le i\le r$. 
Then using the Jacobi identity (\ref{phi-Jacobi-qva}) we obtain
\begin{eqnarray}\label{especial-weakass-n}
(z+1)^{n}Y_{W}(u,(z+1)x_2)Y_{W}(v,x_2)w
=(1+z)^{n}Y_{W}(Y(u,\log (1+z))v,x_2)w.
\end{eqnarray}
Noticing that $(z+1)^{m}Y_{W}(u,(z+1)x_2)Y_{W}(v,x_2)w$ exists for any $m\in \Z$, we have
\begin{eqnarray*}
&&u_{[0]}v_{[0]}w\nonumber\\
&=&\Res_{z}\Res_{x_2}x_2^{-1}(z+1)^{-1}Y_{W}(u,(z+1)x_2)Y_{W}(v,x_2)w\nonumber\\
&=&\Res_{z}\Res_{x_2}x_2^{-1}(z+1)^{-n-1}\left[(z+1)^{n}Y_{W}(u,(z+1)x_2)Y_{W}(v,x_2)w\right]\nonumber\\
&=&\Res_{z}\Res_{x_2}\sum_{j\ge 0}\binom{-n-1}{j}z^{-n-1-j}x_2^{-1}\left[(z+1)^{n}Y_{W}(u,(z+1)x_2)Y_{W}(v,x_2)w\right].
\end{eqnarray*}
We claim
$$\Res_{z}\Res_{x_2}z^{-n-1-j}x_2^{-1}\left[(z+1)^{n}Y_{W}(u,(z+1)x_2)Y_{W}(v,x_2)w\right]=0$$
for $j\ge n+1$. Indeed, for $j\ge n+1$, we have
\begin{eqnarray*}
&&\Res_{z}\Res_{x_2}z^{-n-1-j}x_2^{-1}(z+1)^{n}Y_{W}(u,(z+1)x_2)Y_{W}(v,x_2)w\nonumber\\
&=&\sum_{p,q\in \Z}\Res_{z}\Res_{x_2}z^{-n-1-j}x_2^{-1-p-q}(z+1)^{n-p}u_{[p]}v_{[q]}w\nonumber\\
&=&\sum_{p,q\in \Z}\sum_{i\ge 0}\binom{n-p}{i}\Res_{z}\Res_{x_2}z^{-1-j-p-i}x_2^{-1-p-q}u_{[p]}v_{[q]}w\nonumber\\
&=&\sum_{i\ge 0}\binom{n+j+i}{i}u_{[-j-i]}v_{[j+i]}w,
\end{eqnarray*}
which vanishes because $v_{[i+j]}w=0$ with $i+j\ge n+1$. Then we get
\begin{eqnarray*}
&&u_{[0]}v_{[0]}w\nonumber\\
&=&\Res_{z}\Res_{x_2}\sum_{j=0}^{n}\binom{-n-1}{j}z^{-n-1-j}x_2^{-1}\left[(z+1)^{n}Y_{W}(u,(z+1)x_2)Y_{W}(v,x_2)w\right]
\nonumber\\
&=&\Res_{z}\Res_{x_2}\sum_{j=0}^{n}\binom{-n-1}{j}z^{-n-1-j}x_2^{-1}\left[(1+z)^{n}Y_{W}(Y(u,\log (1+z))v,x_2)w\right]\nonumber\\
&=&(u\tilde{*}_{n}v)_{[0]}w.
\end{eqnarray*}
For proving that $\widetilde{\Omega}_{n}(W)$ is an $\widetilde{A}_{n}(V)$-module, 
it remains to show that $(u\tilde{\circ}_{n}v)_{[0]}w=0$ for $u,v\in V,\ w\in \widetilde{\Omega}_{n}(W)$.
Recall from (\ref{efn-derivative}) that $f_{n}'(x)=\binom{-n-1}{n}(2n+1)\frac{e^{(n+1)x}}{(e^x-1)^{2n+2}}$.
With $u,v\in V,\ w\in \widetilde{\Omega}_{n}(W)$, using the first part we get
\begin{eqnarray*}
&&\binom{-n-1}{n}(2n+1)(u\tilde{\circ}_{n}v)_{[0]}w\nonumber\\
&=&\Res_{x_2}\Res_{x}\binom{-n-1}{n}(2n+1)x_2^{-1}\frac{e^{(n+1)x}}{(e^x-1)^{2n+2}}Y_{W}(Y(u,x)v,x_2)w\\
&=&\Res_{x_2}\Res_{x}x_2^{-1}\left(\frac{d}{dx}f_{n}(x)\right)Y_{W}(Y(u,x)v,x_2)w\\
&=&-\Res_{x_2}\Res_{x}x_2^{-1}f_{n}(x)\frac{\partial}{\partial x}Y_{W}(Y(u,x)v,x_2)w\\
&=&-\Res_{x_2}\Res_{x}x_2^{-1}f_{n}(x)Y_{W}(Y(\D u,x)v,x_2)w\\
&=&-((\D u)*_{n}v)_{[0]}w\\
&=&-(\D u)_{[0]}v_{[0]}w\\
&=&0,
\end{eqnarray*}
noticing that $(\D u)_{[0]}=0$ from the $\D$-derivative property (\ref{ederivative-property}).
Thus $(u\tilde{\circ}_{n}v)_{[0]}w=0$. Therefore,
$\widetilde{\Omega}_{n}(W)$ is an $\widetilde{A}_{n}(V)$-module with $v$ acting as $v_{[0]}$ for $v\in V$.

Let $W=\oplus_{m\in \N}W[m]$ be an $\N$-graded $\phi$-coordinated $V$-module.
Assume $w=\sum w_{m}\in \widetilde{\Omega}_{0}(W)$ with $w_{m}\in W[m]$. It is straightforward to see that 
$w_{m}\in \widetilde{\Omega}_{0}(W)$ for all $m$. This shows that $\widetilde{\Omega}_{0}(W)$ is a graded subspace.
It can be readily seen that $W[0]$ 
is an $\widetilde{A}(V)$-submodule of $\widetilde{\Omega}_{0}(W)$. 
Using weak $\phi$-associativity one can show that for any $w\in W$, the linear span of $v_{[m]}w$ 
for $v\in V,\ m\in \Z$ is a $\phi$-coordinated submodule of $W$ (cf. \cite{DM}, \cite{li-thesis}).
Then for any $\widetilde{A}(V)$-submodule $U$ of $W[0]$,
the linear span of $v_{[m]}U$ for $v\in V,\ m\in \Z$ is a graded $\phi$-coordinated $V$-submodule of $W$ 
with $U$ as the degree-zero subspace.
 It follows that if $W$ is also irreducible, then $\widetilde{\Omega}_{0}(W)=W[0]$ and $W[0]$ must be an irreducible 
 $\widetilde{A}(V)$-module.
\end{proof}

\section{Associating $\N$-graded $\phi$-coordinated $V$-modules to $\widetilde{A}(V)$-modules}

Let $V$ be a weak quantum vertex algebra with a fixed constant $\S$-map. In this section, to $V$ we 
associate an associative algebra $\widetilde{U}(V)$. Under a certain assumption on the structure of $\widetilde{U}(V)$,
for any $\widetilde{A}(V)$-module $U$, we construct an $\N$-graded $\phi$-coordinated $V$-module $L(U)$
with $L(U)[0]\simeq U$ as an $\widetilde{A}(V)$-module. This together with Proposition \ref{pomega}
gives rise to a bijection between the set of isomorphism classes of irreducible $\widetilde{A}(V)$-modules and
 the set of equivalence classes of irreducible $\N$-graded $\phi$-coordinated $V$-modules.
On the other hand, we show that if $V$ is a vertex super-algebra, the very assumption holds true.
In the process of the proof, to $V$ we associate a $\Z$-graded Lie super-algebra ${\mathcal{L}}_{\phi}(V)$, whose  
universal enveloping algebra is shown essentially to be isomorphic to $\widetilde{U}(V)$.
The validity of the assumption follows from  the P-B-W theorem. 

Let $V$ be a weak quantum vertex algebra with a constant $\S$-map, which is fixed throughout this section. 
From now on, we fix a linear operator $\S\in \End _{\C}(V\otimes V)$  such that for $u,v\in V$, the $\S$-Jacobi identity 
(\ref{e2.12}) holds  with $u^{(i)},v^{(i)}$ given by
$$\S(v\otimes u)=\sum_{i=1}^{r}v^{(i)}\otimes u^{(i)}.$$
We also have
$$Y(u,x)v=\sum_{i=1}^{r}e^{x\D}Y(v^{(i)},-x)u^{(i)}.$$

Recall that for a $\phi$-coordinated $V$-module $(W,Y_{W})$, we have
\begin{eqnarray}\label{phi-Jacobi}
&&(xz)^{-1} \delta\left(\frac{x_1-x}{xz}\right)Y_{W}(u, x_1)Y_{W}(v, x) \nonumber\\
&&\hspace{1cm}  -(xz)^{-1}\delta \left(\frac{x-x_1}{-xz}\right)
\sum_{i=1}^{r}Y_{W}(v^{(i)}, x)Y_{W}(u^{(i)}, x_1) \nonumber\\
&&=x_1^{-1}\delta\left(\frac{x(1+z)}{x_1}\right)Y_{W}\left(Y(u, \log(1 + z))v, x\right)
\end{eqnarray}
for $u,v\in V$. From (\ref{phi-Jacobi}) we obtain the following generalized commutator formula
\begin{eqnarray}\label{Lie-commutator-bracket}
&&Y_{W}(u,x_1)Y_{W}(v,x_2)-\sum_{i=1}^{r}Y_{W}(v^{(i)},x_2)Y_{W}(u^{(i)},x_1)\nonumber\\
&=&\Res_{z} x_2x_1^{-1}\delta\left(\frac{x_2(1+z)}{x_1}\right)Y_{W}(Y(u, \log(1 + z))v, x_2)  \nonumber\\
&=&\Res_{x_0}Y_{W}(Y(u,x_0)v,x_2)\delta\left(\frac{x_2e^{x_0}}{x_1}\right).
\end{eqnarray}
Noticing that 
$$\delta\left(\frac{xe^{x_0}}{x_1}\right)=e^{x_0(x_2\frac{\partial}{\partial x_2})}\delta\left(\frac{x_2}{x_1}\right),$$
we get
\begin{eqnarray}\label{ecommutator-phi}
&&Y_{W}(u,x_1)Y_{W}(v,x_2)-\sum_{i=1}^{r}Y_{W}(v^{(i)},x_2)Y_{W}(u^{(i)},x_1)\nonumber\\
&=&\sum_{j\ge 0}\frac{1}{j!}Y_{W}(u_{j}v,x_2)\left(x_2\frac{\partial}{\partial x_2}\right)^{j}\delta\left(\frac{x_2}{x_1}\right).
\end{eqnarray}
In terms of components, we have
\begin{eqnarray}\label{ecommutator-phi-comp}
u_{[m]}v_{[n]}-\sum_{i=1}^{r}v^{(i)}_{[n]}u^{(i)}_{[m]}
=\sum_{j\ge 0}\frac{m^{j}}{j!}\left(u_{j}v\right)_{[m+n]}
\end{eqnarray}
for $u,v\in V,\  m,n\in \Z$.

\bd{dalgebraU}
{\em Let $\widetilde{U}(V)$ be the associative algebra with generators $u[n]$ (linear in $u$) for $u\in V,\ n\in \Z$, 
 subject to relations
\begin{eqnarray}
&& {\bf 1}[n]=\delta_{n,0},\   \   \   \  (\D u)[n]=-nu[n],\label{relation1d}\\
&&u[m]v[n]-\sum_{i=1}^{r}v^{(i)}[n]u^{(i)}[m]=\sum_{j\ge 0}\frac{m^{j}}{j!}(u_{j}v)[m+n].\label{relation-mn}
 \end{eqnarray}}
 \ed

It can be readily seen that defining $\deg u[n]=-n$ for $u\in V,\ n\in \Z$ makes $\widetilde{U}(V)$ a $\Z$-graded algebra.
Then we have the nations of $\Z$-graded $\widetilde{U}(V)$-module and $\N$-graded $\widetilde{U}(V)$-module.
 Set
\begin{eqnarray}
 \widetilde{U}(V)^{0}=\< u[0]\  |\  u\in V\>,  \    \    \   \
 \widetilde{U}(V)^{\pm}=\< u[\mp]\ |\  u\in V,\ n\in \Z_{+}\>,
 \end{eqnarray}
where $\< T\>$ denotes the subalgebra generated by $T$ for a subset $T$. We have  
\begin{eqnarray}\label{etriangular-product}
\widetilde{U}(V)=\widetilde{U}(V)^{+}\widetilde{U}(V)^{0}\widetilde{U}(V)^{-}.
\end{eqnarray}
Notice that $\widetilde{U}(V)^{0}\widetilde{U}(V)^{-}$ is exactly the subalgebra generated by $u[n]$ with $u\in V,\ n\ge 0$.

\bd{dBalgebra}
{\em Denote by $B(V)$ the associative algebra with generators $u[n]$ (linear in $u$) for $u\in V,\ n\in \N$, 
subject to the relations (\ref{relation1d}) and (\ref{relation-mn}) with $m,n\ge 0$.}
\ed

From the definitions of $B(V)$ and $\widetilde{U}(V)$,  there is an algebra homomorphism
\begin{eqnarray}
\pi: \  B(V) \rightarrow \widetilde{U}(V),
\end{eqnarray}
sending $u[n]\in B(V)$ to $u[n]\in  \widetilde{U}(V)$ for $u\in V,\ n\ge 0$ (with abuse of notation), where
$\pi(B(V))=\widetilde{U}(V)^{0}\widetilde{U}(V)^{-}$. 
On the other hand, we have:

\bl{Bfact}
There is an algebra homomorphism from $B(V)$ onto $\widetilde{A}(V)$ such that
\begin{eqnarray}
\psi (u[m])=\delta_{m,0}(u+\widetilde{O}(V))
\end{eqnarray}
for $u\in V,\ m\ge 0$.
\el

\begin{proof}  First, we claim $\D V\subset \widetilde{O}(V)$ and for $u,v\in V$,
\begin{eqnarray}\label{egcommutator}
u*v-\sum_{i=1}^{r}v^{(i)}*u^{(i)}\equiv u_{0}v\   \  \  (\mod\ \D V).
\end{eqnarray}
 Indeed, for $v\in V$, we have
\begin{eqnarray*}
v\widetilde{\circ} {\bf 1}=\Res_{z}\frac{1}{z^2}Y(v,\log (1+z)){\bf 1}
=\Res_{z}\frac{1}{z^2}e^{\D \log (1+z)}v=\Res_{z}\frac{1}{z^2}\log (1+z) \D v
=\D v.
\end{eqnarray*}
This shows that $\D V\subset \widetilde{O}(V)$.
For $u,v\in V$, we have
\begin{eqnarray*}
&&u*v-\sum_{i=1}^{r}v^{(i)}*u^{(i)}\nonumber\\
&=&\Res_{x}\frac{e^{x}}{e^x-1}Y(u,x)v-\Res_{x}\frac{e^{x}}{e^x-1}Y(v^{(i)},x)u^{(i)}\nonumber\\
&=&\Res_{x}\frac{e^{x}}{e^x-1}Y(u,x)v-\Res_{x}\frac{e^{x}}{e^x-1}e^{x\D}Y(u,-x)v\nonumber\\
&=&\Res_{x}\frac{e^{x}}{e^x-1}Y(u,x)v+\Res_{x}\frac{e^{-x}}{e^{-x}-1}e^{-x\D}Y(u,x)v\nonumber\\
&\equiv &\Res_{x}Y(u,x)v\  \   (\mod\; \D V)\nonumber\\
&=&u_{0}v,
 \end{eqnarray*}
 proving (\ref{egcommutator}). Then it is straightforward to show that the desired algebra homomorphism $\psi$
 exists.
 \end{proof}
  
Now, let $U$ be an $\widetilde{A}(V)$-module. Our goal is to construct an $\N$-graded $\phi$-coordinated $V$-module 
$W=\oplus_{n\in \N}W[n]$ with $W[0]\simeq U$ as an $\widetilde{A}(V)$-module.
View $U$ as a $B(V)$-module through the algebra homomorphism $\psi$ obtained in Lemma \ref{Bfact}.
On the other hand, view $\widetilde{U}(V)$ as a $B(V)$-bimodule 
by using the algebra homomorphism $\pi$ from  $B(V)$ to $\widetilde{U}(V)$.
Then form induced module
\begin{eqnarray}
\tilde{M}(U)=\widetilde{U}(V)\otimes _{B(V)}U,
\end{eqnarray}
an $\N$-graded $\widetilde{U}(V)$-module and an $\N$-graded $B(V)$-module with $\deg (1\otimes U)=0$.
Denote by $\theta_{U}$ the linear map
\begin{eqnarray}
\theta_{U}:\  U\rightarrow \tilde{M}(U)[0]\subset \tilde{M}(U);\  \  u\mapsto 1\otimes u,
\end{eqnarray}
which is a $B(V)$-module homomorphism. We hope that $\theta_{U}$ is an embedding, so that
$\tilde{M}(U)$ is an $\N$-graded $\widetilde{U}(V)$-module  with
$\tilde{M}(U)[0]\simeq U$ as a $B(V)$-module.

For this purpose, we make the following assumption (cf. (\ref{etriangular-product})):

{\bf Assumption A:}   $\widetilde{U}(V)\simeq \widetilde{U}(V)^{+}\otimes \widetilde{U}(V)^{0}\widetilde{U}(V)^{-}$
as a vector space.

This assumption is similar to part of the assertion of the celebrated P-B-W theorem. We are not able to show
that it holds for the general case in this paper.
In the following, we prove that Assumption A holds if $V$ is a vertex super-algebra. 

We first associate a Lie super-algebra to a vertex super-algebra $V$ (cf. \cite{ffr}).

\bp{plie-algebra}
Let $V=V_{\bar{0}}\oplus V_{\bar{1}}$ be a vertex super-algebra. Define a bilinear operation $[\cdot,\cdot]_{\phi}$ on vector space $V\otimes \C[t,t^{-1}]$ by
\begin{eqnarray}\label{ephi-commutator}
[u\otimes t^{m},v\otimes t^{n}]_{\phi}=\sum_{j\ge 0}\frac{m^{j}}{j!}u_{j}v\otimes t^{m+n}
\end{eqnarray}
for $u,v\in V,\ m,n\in \Z$. Set 
$$I=\span\{ \D v\otimes t^{n}+n(v\otimes t^{n})\ |\ v\in V,\ n\in \Z\}.$$
Then $I$ is a (two-sided) ideal of the non-associative algebra $(V\otimes \C[t,t^{-1}], [\cdot,\cdot]_{\phi})$
and  the quotient algebra is a Lie super-algebra,
which is denoted by ${\mathcal{L}}_{\phi}(V)$.
\ep

\begin{proof} Consider the Laurent polynomial algebra $\C[t,t^{-1}]$. Set $D=t\frac{d}{dt}$, a derivation on $\C[t,t^{-1}]$.
By a result of Borcherds, $\C[t,t^{-1}]$ becomes a vertex algebra 
with $1$ as the vacuum vector and with
\begin{eqnarray}
Y(f(t),x)g(t)=\left(e^{xD}f(t)\right)g(t)=f(te^{x})g(t)
\end{eqnarray}
for $f(t),g(t)\in \C[t,t^{-1}]$.
For convenience, we denote this vertex super-algebra by $\C[t,t^{-1}]_{\phi}$. 
Furthermore, with vertex super-algebra $V$, we have a tensor product 
vertex super-algebra $V\otimes \C[t,t^{-1}]_{\phi}$, where
\begin{eqnarray}
Y(u\otimes f(t),x)(v\otimes g(t))&=&Y(u,x)v\otimes Y(f(t),x)g(t)\nonumber\\
&=&Y(u,x)v\otimes f(te^{x})g(t)
\end{eqnarray}
for $u,v\in V,\ f(t),g(t)\in \C[t,t^{-1}]$ and where $V_{\bar{i}}\otimes \C[t,t^{-1}]$ for $i=0,1$ are the even part and the odd part, respectively.
A result of Borcherds states that  for each vertex super-algebra $K$, $\D K$ is a two-sided ideal of
the non-associative algebra $(K,*)$, where $u*v=u_{0}v$ for $u,v\in K$, and the quotient  $K/\D K$ is a super-algebra.
Notice that the $\D$-operator of $V\otimes \C[t,t^{-1}]_{\phi}$ is
$$\hat{\D}=\D_{V}\otimes 1+1\otimes t\frac{d}{dt}.$$
Taking $K=V\otimes \C[t,t^{-1}]_{\phi}$, we obtain a Lie super-algebra structure
on the quotient space $(V\otimes \C[t,t^{-1}])/\hat{\D}(V\otimes \C[t,t^{-1}])$, where
$$(u\otimes t^{m})_{0}(v\otimes t^{n})=\sum_{j\ge 0}\frac{m^{j}}{j!}u_{j}v\otimes t^{m+n}
=[u\otimes t^{m},v\otimes t^{n}]_{\phi}$$
for $u,v\in V,\ m,n\in \Z$, as 
$$(u\otimes t^{m})_{0}=\Res_{x}Y(u\otimes t^{m},x)=\Res_{x}(Y(u,x)\otimes t^{m}e^{mx})
=\sum_{j\ge 0}\frac{m^{j}}{j!}(u_{j}\otimes t^{m}).$$
Then the assertions follow immediately.
\end{proof}

Assume that $V$ is a vertex super-algebra. 
It can be readily seen that ${\mathcal{L}}_{\phi}(V)$ is a $\Z$-graded Lie super-algebra with 
\begin{eqnarray}
\deg v(n)=-n\   \   \   \mbox{  for }v\in V,\ n\in \Z,
\end{eqnarray}
where $v(n)$ denotes the image of $v\otimes t^{n}$ in ${\mathcal{L}}_{\phi}(V)$.
For $n\in \Z$, denote by ${\mathcal{L}}_{\phi}(V)(n)$ the
homogeneous subspace of degree $n$. Set
\begin{eqnarray}
{\mathcal{L}}_{\phi}(V)_{\pm}=\oplus_{n\in \Z_{+}}{\mathcal{L}}_{\phi}(V)(\mp n).
\end{eqnarray}

Recall that ${\bf 1}[n]=\delta_{n,0}$ in $\widetilde{U}(V)$ for $n\in \Z$.
On the other hand, we have  that ${\bf 1}(n)=0$ for $n\ne 0$ and ${\bf 1}(0)$ is central in ${\mathcal{L}}_{\phi}(V)$,
noticing that $\D {\bf 1}=0$. It follows straightforwardly that
\begin{eqnarray}
\widetilde{U}(V)\simeq U({\mathcal{L}}_{\phi}(V))/I,\   \  \  \   
B(V)\simeq U({\mathcal{L}}_{\phi}(V)_{-}+{\mathcal{L}}_{\phi}(V)(0))/I',
\end{eqnarray}
where $I=U({\mathcal{L}}_{\phi}(V))({\bf 1}(0)-1)$,  $I'=U({\mathcal{L}}_{\phi}(V)_{-}+{\mathcal{L}}_{\phi}(V)(0))({\bf 1}(0)-1)$, 
which are ideals of the corresponding algebras. 
Furthermore, by the P-B-W theorem we have
\begin{eqnarray}
\widetilde{U}(V)&\simeq& U({\mathcal{L}}_{\phi}(V)_{+})\otimes U({\mathcal{L}}_{\phi}(V)(0))/I_0
\otimes U({\mathcal{L}}_{\phi}(V)_{-}),\\
B(V)&\simeq &U({\mathcal{L}}_{\phi}(V)(0))/I_0 \otimes U({\mathcal{L}}_{\phi}(V)_{-}),
\end{eqnarray}
where $I_0=U({\mathcal{L}}_{\phi}(V)(0))({\bf 1}(0)-1)$. It follows that
$$\widetilde{U}(V)^{\pm}\simeq U({\mathcal{L}}_{\phi}(V)_{\pm}),\    \   \   \   
\widetilde{U}(V)^{0}\simeq U({\mathcal{L}}_{\phi}(V)(0))/I_0.$$
As an immediate consequence have:

\bc{cuniversal}
Assume that $V$ is a vertex super-algebra. Then 
$$ \widetilde{U}(V)=\widetilde{U}(V)^{+}\otimes \widetilde{U}(V)^{0}\widetilde{U}(V)^{-}.$$
\ec

For the rest of this section, we assume that {\bf Assumption A} holds. This particularly implies that
the linear map $\theta_{U}: U\rightarrow  \tilde{M}(U)[0]$ is a $B(V)$-module isomorphism.
Then we identify $\tilde{M}(U)[0]$ with $U$ this way. 

Let $P_{0}$ denote the projection map of $\tilde{M}(U)$ onto $U$ with respect to the $\N$-graded decomposition.
For $u^{*}\in U^{*}, \ w\in \tilde{M}(U)$, we define
\begin{eqnarray}
\< u^{*},w\>=\<u^{*},P_{0}(w)\>.
\end{eqnarray}
Set 
$$J(U)=\left\{ w\in \tilde{M}(U)\ |\ \<u^{*},aw\>=0\  \  \mbox{ for all }u^{*}\in U^{*},\ a\in \widetilde{U}(V)\right\}.$$
Then $J(U)$ is a graded $\widetilde{U}(V)$-submodule of $\tilde{M}(U)$. Furthermore, we set 
\begin{eqnarray}
L(U)=\tilde{M}(U)/J(U),
\end{eqnarray}
which is an $\N$-graded  $\widetilde{U}(V)$-module. 

The following is the main result: 

\bt{train}
 $L(U)$ is an $\N$-graded $\phi$-coordinated $V$-module with $U$ as the degree zero subspace.
 Furthermore, if $U$ is an irreducible $\widetilde{A}(V)$-module, then $L(U)$ is an irreducible $\N$-graded 
 $\phi$-coordinated $V$-module.
 \et
 
\begin{proof}  It is clear that $J(U)\cap U=0$. Thus $L(U)$ is an $\N$-graded $\widetilde{U}(V)$-module with $U$ as 
the degree zero subspace. It remains to show that $L(U)$ is a $\phi$-coordinated $V$-module. 
We need to establish the $\phi$-coordinated-module Jacobi identity.

{\bf Claim 1:} For $u,v\in V,\ w^{*}\in U^{*},\ w\in U$, 
\begin{eqnarray}\label{eYWuYWv-1}
&&(x_1-x_2)^{-1}\<w^{*},Y_{W}(u,x_1)Y_{W}(v,x_2)w\>\nonumber\\
&=&\Res_{x_0}\frac{e^{x_0}}{(e^{x_0}-1)(x_1-x_2e^{x_0})}\<w^{*},Y_{W}(Y(u,x_0)v,x_2)w\>\nonumber\\
&=&\Res_{z}z^{-1}(x_1-x_2(1+z))^{-1}\<w^{*},Y_{W}(Y(u,\log (1+z))v,x_2)w\>.
\end{eqnarray}

 Noticing that $P_{0}(u_{[m]}v_{[n]}w)=0$ for $m,n\in \Z$ with $m+n\ne 0$
 and that $v_{[n]}w=0$ for $n\ge 1$, we have
\begin{eqnarray*}
&&\<w^{*},Y_{W}(u,x_1)Y_{W}(v,x_2)w\>\\
&=&\sum_{m\in \Z}\<w^{*},u_{[m]}v_{[-m]}w\>x_1^{-m}x_2^{m}\\
&=&\sum_{m\ge 0}\<w^{*},u_{[m]}v_{[-m]}w\>x_1^{-m}x_2^{m}\\
&=&\<w^{*},u_{[0]}v_{[0]}w\>+\sum_{m\ge 1}\<w^{*},u_{[m]}v_{[-m]}w\>x_1^{-m}x_2^{m}\\
&=&\<w^{*},u_{[0]}v_{[0]}w\>+\sum_{m\ge 1,j\ge 0}\frac{m^{j}}{j!}\<w^{*},(u_{j}v)_{[0]}w\>x_1^{-m}x_2^{m}.
\end{eqnarray*}
 Furthermore, we have
\begin{eqnarray*}
&&\sum_{m\ge 1,j\ge 0}\frac{m^{j}}{j!}\<w^{*},(u_{j}v)_{[0]}w\>x_1^{-m}x_2^{m}\\
&=&\sum_{m\ge 1,j\ge 0}\frac{1}{j!}\left(x_2\frac{\partial}{\partial x_2}\right)^{j}(x_2/x_1)^{m}\<w^{*},(u_{j}v)_{[0]}w\>\\
&=&\sum_{j\ge 0}\frac{1}{j!}\left(x_2\frac{\partial}{\partial x_2}\right)^{j}\frac{x_2}{x_1-x_2}\<w^{*},(u_{j}v)_{[0]}w\>\\
&=&\sum_{j\ge 0}\left(\frac{1}{j!}\left(x_2\frac{\partial}{\partial x_2}\right)^{j}\frac{x_2}{x_1-x_2}\right)
\<w^{*},Y_{W}(u_{j}v,x_2)w\>\\
&=&\Res_{x_0}\left(e^{x_0x_2\frac{\partial}{\partial x_2}}\frac{x_2}{x_1-x_2}\right)\<w^{*},Y_{W}(Y(u,x_0)v,x_2)w\>\\
&=&\Res_{x_0}\frac{x_2e^{x_0}}{x_1-x_2e^{x_0}}\<w^{*},Y_{W}(Y(u,x_0)v,x_2)w\>,
\end{eqnarray*}
noticing that $\<w^{*},Y_{W}(u_{j}v,x_2)w\>=\<w^{*},(u_{j}v)_{[0]}w\>$. 
On the other hand, as $U$ is an $\widetilde{A}(V)$-module we have
\begin{eqnarray*}
\<w^{*},u_{[0]}v_{[0]}w\>=\Res_{x_0}\frac{e^{x_0}}{e^{x_0}-1}\<w^{*},Y_{W}(Y(u,x_0)v,x_2)w\>.
\end{eqnarray*}
Then we obtain
\begin{eqnarray*}
&&\<w^{*},Y_{W}(u,x_1)Y_{W}(v,x_2)w\>\nonumber\\
&=&\Res_{x_0}\frac{e^{x_0}}{e^{x_0}-1}\<w^{*},Y_{W}(Y(u,x_0)v,x_2)w\>\nonumber\\
&&+\Res_{x_0}\frac{x_2e^{x_0}}{x_1-x_2e^{x_0}}\<w^{*},Y_{W}(Y(u,x_0)v,x_2)w\>\nonumber\\
&=&\Res_{x_0}\frac{ (x_1-x_2)e^{x_0}}{(e^{x_0}-1)(x_1-x_2e^{x_0})}\<w^{*},Y_{W}(Y(u,x_0)v,x_2)w\>,
\end{eqnarray*}
which gives (\ref{eYWuYWv-1}). 

{\bf Claim 2:} For $u,v\in V,\ w^{*}\in U^{*},\  w\in U$,
\begin{eqnarray}\label{cross-prod-1}
&&(x_1-x_2)^{-1}\<w^{*},Y_{W}(u,x_1)Y_{W}(v,x_2)w\>\nonumber\\
&&\hspace{0.5cm}-(-x_2+x_1)^{-1}\sum_{i=1}^{r}\<w^{*},Y_{W}(v^{(i)},x_2)Y_{W}(u^{(i)},x_1)w\>\nonumber\\
&=&\Res_{z}z^{-1}x_1^{-1}\delta\left(\frac{x_2(1+z)}{x_1}\right)\<w^{*},Y_{W}(Y(u,\log(1+z))v,x_2)w\>, 
\end{eqnarray}
where $u^{(i)},v^{(i)}\in V$, $i=1,\dots, r$ such that 
$$Y(u,x)v=\sum_{i=1}^{r}e^{x\D}Y(v^{(i)},-x)u^{(i)}.$$
From Claim 1, we have
\begin{eqnarray}
&&(x_2-x_1)^{-1}\sum_{i=1}^{r}\<w^{*},Y_{W}(v^{(i)},x_2)Y_{W}(u^{(i)},x_1)w\>\nonumber\\
&=&-\Res_{x_0}\frac{e^{-x_0}}{(e^{-x_0}-1)(x_2-x_1e^{-x_0})}\sum_{i=1}^{r}\<w^{*},Y_{W}(Y(v^{(i)},-x_0)u^{(i)},x_1)w\>\nonumber\\
&=&\Res_{x_0}\frac{e^{x_0}}{(e^{x_0}-1)(x_2e^{x_0}-x_1)}\sum_{i=1}^{r}
\<w^{*},Y_{W}(Y(v^{(i)},-x_0)u^{(i)},x_1)w\>\nonumber\\
&=&\Res_{x_0}\frac{e^{x_0}}{(e^{x_0}-1)(x_2e^{x_0}-x_1)}\<w^{*},Y_{W}(e^{-x_0\D}Y(u,x_0)v,x_1)w\>\nonumber\\
&=&\Res_{x_0}\frac{e^{x_0}}{(e^{x_0}-1)(x_2e^{x_0}-x_1)}\<w^{*},Y_{W}(Y(u,x_0)v,x_1e^{-x_0})w\>\nonumber\\
&=&\Res_{x_0}\frac{e^{x_0}}{(e^{x_0}-1)(x_2e^{x_0}-x_1)}\<w^{*},Y_{W}(Y(u,x_0)v,x_2)w\>,
\end{eqnarray}
noticing that $\<w^{*},Y_{W}(Y(u,x_0)v,x_2)w\>$ does not depend on $x_2$.
Combining Claim 1 with this  we obtain (\ref{cross-prod-1}) as
\begin{eqnarray*}
&&(x_1-x_2)^{-1}\<w^{*},Y_{W}(u,x_1)Y_{W}(v,x_2)w\>\nonumber\\
&&\hspace{0.5cm}-(-x_2+x_1)^{-1}\sum_{i=1}^{r}\<w^{*},Y_{W}(v^{(i)},x_2)Y_{W}(u^{(i)},x_1)w\>\nonumber\\
&=&\Res_{x_0}\frac{e^{x_0}}{e^{x_0}-1}\left((x_1-x_2e^{x_0})^{-1}-(-x_2e^{x_0}+x_1)^{-1}\right)
\<w^{*},Y_{W}(Y(u,x_0)v,x_2)w\>
\nonumber\\
&=&\Res_{x_0}\frac{e^{x_0}}{e^{x_0}-1}x_1^{-1}\delta\left(\frac{x_2e^{x_0}}{x_1}\right)\<w^{*},Y_{W}(Y(u,x_0)v,x_2)w\>
\nonumber\\
&=&\Res_{z}z^{-1}x_1^{-1}\delta\left(\frac{x_2(1+z)}{x_1}\right)\<w^{*},Y_{W}(Y(u,\log(1+z))v,x_2)w\>.
\end{eqnarray*}

{\bf Claim 3:} For $u,v\in V,\ w^{*}\in U^{*},\  w\in U$,
\begin{eqnarray}\label{eclaim3}
&&(zx_2)^{-1}\delta\left(\frac{x_1-x_2}{zx_2}\right)\<w^{*},Y_{W}(u,x_1)Y_{W}(v,x_2)w\>\nonumber\\
&&\hspace{0.5cm}-(zx_2)^{-1}\delta\left(\frac{x_2-x_1}{-zx_2}\right)
\sum_{i=1}^{r}\<w^{*},Y_{W}(v^{(i)},x_2)Y_{W}(u^{(i)},x_1)w\>\nonumber\\
&=&x_1^{-1}\delta\left(\frac{x_2(1+z)}{x_1}\right)\<w^{*},Y_{W}(Y(u,\log(1+z))v,x_2)w\>, 
\end{eqnarray}
where $u^{(i)},v^{(i)}\in V$, $i=1,\dots, r$, are given as before. 

Notice that (\ref{eclaim3}) is equivalent to
\begin{eqnarray}\label{en-jacobi}
&&(x_1-x_2)^{n}\<w^{*},Y_{W}(u,x_1)Y_{W}(v,x_2)w\>\nonumber\\
&&\   \   \   -(-x_2+x_1)^{n}\sum_{i=1}^{r}\<w^{*},Y_{W}(v^{(i)},x_2)Y_{W}(u^{(i)},x_1)w\>\nonumber\\
&=&\Res_{z}x_2(zx_2)^{n}x_1^{-1}\delta\left(\frac{x_2(1+z)}{x_1}\right)\<w^{*},Y_{W}(Y(u,\log(1+z))v,x_2)w\>
\end{eqnarray}
for every $n\in \Z$. For $n\ge 0$, this follows from the commutation relation (\ref{Lie-commutator-bracket})
as $\tilde{M}(U)$ is a $\widetilde{U}(V)$-module.
Claim 2 asserts that it is true for $n=-1$. We then use induction to show that it is true for all negative integers $n$.
The induction goes as follows: Suppose that $n$ is a negative integer such that  (\ref{en-jacobi}) holds.
Applying $x_1\frac{\partial}{\partial x_1}$ to (\ref{en-jacobi}), by using the fact that $(x\frac{d}{dx})Y_{W}(u,x)=Y_{W}(\D u,x)$ we get
\begin{eqnarray}
&&nx_1(x_1-x_2)^{n-1}\<w^{*},Y_{W}(u,x_1)Y_{W}(v,x_2)w\>\nonumber\\
&&\   \   \   \   +(x_1-x_2)^{n}\<w^{*},Y_{W}(\D u,x_1)Y_{W}(v,x_2)w\>\nonumber\\
&&\   -\sum_{i=1}^{r}nx_1(-x_2+x_1)^{n-1}\<w^{*},Y_{W}(v^{(i)},x_2)Y_{W}(u^{(i)},x_1)w\>\nonumber\\
&&-\sum_{i=1}^{r}(-x_2+x_1)^{n}\<w^{*},Y_{W}(v^{(i)},x_2)Y_{W}(\D u^{(i)},x_1)w\>
\nonumber\\
&=&x_1\frac{\partial}{\partial x_1}\Res_{z}x_2(zx_2)^{n}x_1^{-1}\delta\left(\frac{x_2(1+z)}{x_1}\right)
\<w^{*},Y_{W}(Y(u,\log(1+z))v,x_2)w\>\   \   \   \nonumber\\
&=&-\Res_{z}(zx_2)^{n}\left(\frac{\partial}{\partial z}\delta\left(\frac{x_2(1+z)}{x_1}\right)\right)\<w^{*},Y_{W}(Y(u,\log(1+z))v,x_2)w\>\   \   \   \nonumber\\
&=&\Res_{z}nx_2(zx_2)^{n-1}\delta\left(\frac{x_2(1+z)}{x_1}\right)
\<w^{*},Y_{W}(Y(u,\log(1+z))v,x_2)w\>\   \   \   \nonumber\\
&&+\Res_{z}(zx_2)^{n}\delta\left(\frac{x_2(1+z)}{x_1}\right)\frac{\partial}{\partial z}\<w^{*},Y_{W}(Y(u,\log(1+z))v,x_2)w\> \nonumber\\
&=&\Res_{z}nx_2(zx_2)^{n-1}\delta\left(\frac{x_2(1+z)}{x_1}\right)
\<w^{*},Y_{W}(Y(u,\log(1+z))v,x_2)w\>\   \   \   \nonumber\\
&&+\Res_{z}(zx_2)^{n}\delta\left(\frac{x_2(1+z)}{x_1}\right) \frac{1}{1+z}\<w^{*},Y_{W}(Y(\D u,\log(1+z))v,x_2)w\> \nonumber\\
&=&\Res_{z}nx_2(zx_2)^{n-1}\delta\left(\frac{x_2(1+z)}{x_1}\right)
\<w^{*},Y_{W}(Y(u,\log(1+z))v,x_2)w\>\   \   \   \nonumber\\
&&+\Res_{z}x_2(zx_2)^{n}x_1^{-1}\delta\left(\frac{x_2(1+z)}{x_1}\right)\<w^{*},Y_{W}(Y(\D u,\log(1+z))v,x_2)w\>.
\end{eqnarray}
Then it can be readily seen that the case for $n-1$ is also true.

{\bf Claim 4:}  (\ref{eclaim3}) holds for $u,v\in V,\ w^{*}\in U^{*},$ and for all $w\in \tilde{M}(U)$.

Notice that $\tilde{M}(U)$ is generated from $U$ by operators $a_{[m]}$ for $a\in V,\ m<0$.
Suppose $w\in \tilde{M}(U)$ such that (\ref{eclaim3}) holds for $u,v\in V,\ w^{*}\in U^{*}$. Let $a\in V$, $m\in \Z$ with $m<0$.
We have
\begin{eqnarray*}
&&\<w^{*},Y_{W}(u,x_1)Y_{W}(v,x_2)a_{[m]}w\>\nonumber\\
&=&\<w^{*},a_{[m]}Y_{W}(u,x_1)Y_{W}(v,x_2)w\>
-\sum_{j\ge 0}\frac{m^{j}}{j!}x_1^{m}\<w^{*},Y_{W}(a_{j}u,x_1)Y_{W}(v,x_2)w\>\nonumber\\
&&-\sum_{j\ge 0}\frac{m^{j}}{j!}x_2^{m}\<w^{*},Y_{W}(u,x_1)Y_{W}(a_{j}v,x_2)w\>\nonumber\\
&=&-\sum_{j\ge 0}\frac{m^{j}}{j!}x_1^{m}\<w^{*},Y_{W}(a_{j}u,x_1)Y_{W}(v,x_2)w\>\nonumber\\
&&-\sum_{j\ge 0}\frac{m^{j}}{j!}x_2^{m}\<w^{*},Y_{W}(u,x_1)Y_{W}(a_{j}v,x_2)w\>,
\end{eqnarray*}
\begin{eqnarray*}
&&\<w^{*},Y_{W}(Y(u,\log (1+z))v,x_2)a_{[m]}w\>\nonumber\\
&=&\<w^{*},a_{[m]}Y_{W}(Y(u,\log (1+z))v,x_2)w\>\nonumber\\
&&-\sum_{j\ge 0}\frac{m^{j}}{j!}x_2^{m}\<w^{*},Y_{W}(a_{j}Y(u,\log (1+z))v,x_2)w\>\nonumber\\
&=&-\sum_{j\ge 0}\frac{m^{j}}{j!}x_2^{m}\<w^{*},Y_{W}(Y(u,\log (1+z))a_jv,x_2)w\>\nonumber\\
&&-\sum_{i,j\ge 0}\frac{m^{j}}{j!}x_2^{m}{j\choose i}(\log (1+z))^{j-i}\<w^{*},Y_{W}(Y(a_i u,\log (1+z))v,x_2)w\>\nonumber\\
&=&-\sum_{j\ge 0}\frac{m^{j}}{j!}x_2^{m}\<w^{*},Y_{W}(Y(u,\log (1+z))a_jv,x_2)w\>\nonumber\\
&&-\sum_{i\ge 0}\frac{m^{i}}{i!}x_2^{m}(1+z)^{m}\<w^{*},Y_{W}(Y(a_i u,\log (1+z))v,x_2)w\>,
\end{eqnarray*}
noticing that
\begin{eqnarray*}
\sum_{j\ge 0}\frac{m^{j}}{j!}x_2^{m}{j\choose i}(\log (1+z))^{j-i}
&=&\sum_{j\ge i}\frac{m^{i}}{i!}x_2^{m}\frac{1}{(j-i)!}(m\log (1+z))^{j-i}\\
&=&\frac{m^{i}}{i!}x_2^{m}(1+z)^{m}.
\end{eqnarray*}
It can be readily seen that (\ref{eclaim3}) holds with $a_{[m]}w$ in place of $w$.

{\bf Claim 5:} For $u,v\in V$, $w^{*}\in U^{*},\ w\in \tilde{M}(U)$ and for all 
$X\in \widetilde{U}(V)$,
\begin{eqnarray}\label{eclaim5}
&&(zx_2)^{-1}\delta\left(\frac{x_1-x_2}{zx_2}\right)\<w^{*},X\cdot Y_{W}(u,x_1)Y_{W}(v,x_2)w\>\nonumber\\
&&\hspace{0.5cm}-(zx_2)^{-1}\delta\left(\frac{x_2-x_1}{-zx_2}\right)
\sum_{i=1}^{r}\<w^{*},X\cdot Y_{W}(v^{(i)},x_2)Y_{W}(u^{(i)},x_1)w\>\nonumber\\
&=&x_1^{-1}\delta\left(\frac{x_2(1+z)}{x_1}\right)\<w^{*},X\cdot Y_{W}(Y(u,\log(1+z))v,x_2)w\>,
\end{eqnarray}
where $u^{(i)},v^{(i)}\in V$, $i=1,\dots, r$, are given as before. 

Let $X\in \widetilde{U}(V)$ be such that (\ref{eclaim5}) holds.
Let $a\in V,\ m\in \Z$. Consider
$$\<w^{*}, Xa_{[m]}Y_{W}(u,x_1)Y_{W}(v,x_2)w\>,$$
$$\<w^{*}, Xa_{[m]}Y_{W}(v,x_2)Y_{W}(u,x_1)w\>,$$
$$\<w^{*}, Xa_{[m]}Y_{W}(Y(u,\log(1+z))v,x_2)w\>.$$
From the proof of Claim 4, we see that (\ref{eclaim5}) holds with $Xa_{[m]}$ in place of $X$.
It then follows from induction.

It is straightforward to show that if $U$ is an irreducible $\widetilde{A}(V)$-module, 
then $L(U)$ is an irreducible $\N$-graded $\phi$-coordinated $V$-module.
\end{proof}

Recall that every $\phi$-coordinated $V$-module is naturally a $\widetilde{U}(V)$-module.
Let $J$ be the $\widetilde{U}(V)$-submodule of of $\tilde{M}(U)$,
generated by the relation
\begin{eqnarray}\label{edef-J}
&&(xz)^{-1} \delta\left(\frac{x_1-x}{xz}\right)\widetilde{Y}(u,x_1)\widetilde{Y}(v,x)w\nonumber\\
&&\hspace{1cm}-(xz)^{-1}\delta \left(\frac{x-x_1}{-xz}\right)
\sum_{i=1}^{r}\widetilde{Y}(v^{(i)},x)\widetilde{Y}(u^{(i)},x_1)w \nonumber\\
&&
=x_1^{-1}\delta\left(\frac{x(1+z)}{x_1}\right)\widetilde{Y}\left(Y(u, \log(1 + z))v, x\right)w
\end{eqnarray}
for $u,v\in V,\ w\in U\subset \tilde{M}(U)$, where $\widetilde{Y}(u,x)=\sum_{n\in \Z}u[n]x^{-n}\in \widetilde{U}(V)[[x,x^{-1}]]$. 
It can be readily seen that $J$ is a graded $\widetilde{U}(V)$-submodule of $\tilde{M}(U)$.
Then we define
\begin{eqnarray}
M(U)=\tilde{M}(U)/J,
\end{eqnarray}
an $\N$-graded $\widetilde{U}(V)$-module. 

\bp{pM(U)}
Let $U$ be an $\widetilde{A}(V)$-module. Then $M(U)$ is an $\N$-graded $\phi$-coordinated $V$-module
with $M(U)[0]=U$. Furthermore,  for any $\phi$-coordinated $V$-module $W$ and for any 
$A_{\phi}(V)$-module homomorphism $\pi_{0}: U\rightarrow \Omega(W)$, there exists a homomorphism $\pi:
M(U)\rightarrow W$ of $\phi$-coordinated $V$-modules, extending $\pi_0$ uniquely.
\ep

\begin{proof} Recall that $L(U)$ is an $\N$-graded $\phi$-coordinated $V$-module with $L(U)[0]=U$.
From the construction, $L(U)$ as a $\widetilde{U}(V)$-module  is a homomorphism image of $\tilde{M}(U)$.
Then it follows that $M(U)[0]=U$.
From definition, (\ref{edef-J}) holds for $u,v\in V$ 
and for $w\in U\subset M(U)$.
It follows from the proof for Claim 4 in the proof of Theorem \ref{train} that (\ref{edef-J}) holds for $u,v\in V$ 
and for {\em all} $w\in M(U)$.
Then $M(U)$ is an $\N$-graded $\phi$-coordinated $V$-module.
The furthermore-assertion is clear. 
\end{proof}

\section{Rationality, regularity, and fusion rules for vertex operator algebras}
In this section, we study the dependence of rationality, regularity, and fusion  rule on the choice of a conformal vector 
for vertex operator algebras. As the main result, we show that all of these are independent of 
the choice of a conformal vector.

First, define a {\em conformal vertex algebra} (c.f. \cite{hlz}) to be a vertex algebra $V$ equipped with 
 a {\em conformal vector} $\omega$ in the sense that $\omega$ is a vector of $V$ such that
\begin{eqnarray}
[L(m),L(n)]=(m-n)L(m+n)+\frac{1}{12}(m^3-m)\delta_{m+n,0}c
\end{eqnarray}
for $m,n\in \Z$, where $Y(\omega,x)=\sum_{n\in \Z}L(n)x^{-n-2}$,  $c\in \C$, and such that
$L(-1)=\D$, $L(0)$ is semi-simple, and $L(n)$ for $n\ge 1$ are locally nilpotent.
(Recall that $\D$ is the linear operator on $V$ defined by $\D v=v_{-2}{\bf 1}$ for $v\in V$.) Thus
\begin{eqnarray}
[L(-1), Y(v,x)]=Y(L(-1) v,x)=\frac{d}{dx}Y(v,x) \   \   \   \mbox{ for }v\in V.
\end{eqnarray}
The eigenvalues of $L(0)$ on $V$ are called the {\em conformal weights}.
To emphasize the dependence on the conformal vector, we shall also denote a conformal vertex algebra 
by a pair $(V,\omega)$, 

Note that here we do {\em not} assume that a conformal vertex algebra has only integer conformal weights. 
On the other hand, by a {\em conformal vector $\omega$ with integer conformal weights}, 
we mean that all conformal weights are integers.

Recall (see \cite{flm}, \cite{fhl}) that a {\em vertex operator algebra} is a vertex algebra $V$ equipped with a conformal 
vector $\omega$ with integer conformal weights such that
\begin{eqnarray*}
V_{(n)}=0\    \   \mbox{ for $n$ sufficiently negative}
\end{eqnarray*}
and $\dim V_{(n)}<\infty$ for all $n\in \Z$, where $V_{(n)}=\{ v\in V\ |\ L(0)v=nv\}$.

Let $V$ be a vertex operator algebra. A module for $V$ viewed as a vertex algebra is called a {\em weak} $V$-module.
A {\em $V$-module} is a weak $V$-module $W$ such that
$W=\oplus_{\lambda\in \C}W_{(\lambda)}$, 
$\dim W_{(\lambda)}<\infty$ for all $\lambda\in \C$, 
and for any $\lambda\in \C$, $W_{(\lambda+n)}=0$ for sufficiently negative integers $n$,
where $W_{(\lambda)}=\{ w\in W\ |\ L(0)w=\lambda w\}$.

A {\em $\Z$-graded weak $V$-module} is a weak $V$-module $W$ together with a $\Z$-grading 
$W=\oplus_{n\in \Z}W(n)$ such that
\begin{eqnarray}
v_{m}W(n)\subset W(n+\wt v-m-1)
\end{eqnarray}
for homogeneous vector $v\in V$ and for $m,n\in \Z$. A $\Z$-graded weak $V$-module $W=\oplus_{n\in \Z}W(n)$
with $W(n)=0$ for $n<0$ is called an $\N$-graded weak $V$-module.
A weak $V$-module $W$ is said to be {\em admissible} 
if there exists an $\N$-grading with which $W$ becomes an $\N$-graded weak $V$-module. 
 An admissible weak $V$-module is also called an
{\em $\N$-gradable weak $V$-module}  in literature.
Two $\Z$-graded weak $V$-modules $W_i=\oplus_{n\in \Z}W_{i}(n)$ for $i=1,2$ are said to be {\em equivalent} if 
there exists an isomorphism $\psi: W_1\rightarrow W_2$ of weak $V$-modules such that
$\psi(W_1(n))=W_{2}(n+k)$ for $n\in \N$, where $k$ is a fixed integer. Then any nonzero admissible weak $V$-module
is equivalent to an $\N$-graded weak $V$-module $W=\oplus_{n\in \N}W(n)$ with $W(0)\ne 0$.

The following result was due to Zhu (see \cite{zhu1}):

\bl{lzhu}
Let $V$ be a vertex operator algebra and let 
$W=\oplus_{n\in \N}W(n)$ be an irreducible $\N$-graded weak $V$-module with $W(0)\ne 0$.
Then $W=\oplus_{n\in \N}W_{(\lambda+n)}$ for some $\lambda\in \C$, 
where $W_{(\lambda+n)}=\{ w\in W\ |\ L(0)w=(\lambda+n) w\}$.
\el

A vertex operator algebra $V$ is said to be {\em rational}  (see \cite{dlm-twisted}; cf. \cite{zhu1})  
 if every admissible weak module (namely $\N$-gradable weak $V$-module) is completely reducible. 
$V$ is said to be {\em regular} (see \cite{dlm-reg}) if every weak $V$-module is a direct sum of irreducible
(ordinary) $V$-modules. 

A vertex algebra $V$ is said to be {\em $C_2$-cofinite} 
 if $\dim V/C_2(V)<\infty$, where $C_{2}(V)$ is the linear span of $u_{-2}v$ for $u,v\in V$ (see \cite{zhu1}).
It was proved in \cite{li-some} that every regular vertex operator algebra $V$ is $C_{2}$-cofinite (and rational from definition).
Abe, Buhl and Dong  in \cite{abd} proved that the converse also holds, i.e.,  if a vertex operator algebra $V$ is rational and $C_{2}$-cofinite, then it is regular.
On the other hand,  Dong and Yu (see \cite{dy}) proved that a vertex operator algebra $V$ is regular if 
 every $\Z$-graded weak $V$-module  is completely reducible.

Let $V$ be a vertex operator algebra and let $n$ be a nonnegative integer. 
Following \cite{dlm-anv} (page 63), define a binary operation $*_{n}$ on $V$ by  
\begin{eqnarray}
u*_{n}v=\sum_{k=0}^{n}(-1)^{k}\binom{n+k}{n}\Res_{x}Y(u,x)v \frac{(1+x)^{\wt u+n}}{x^{n+k+1}}
\end{eqnarray}
for $u,v\in V$ with $u$ homogeneous. Denote by
$O_{n}(V)$ the subspace of $V$, linearly spanned by $(L(-1)+L(0))V$ and by
\begin{eqnarray}
u\circ_{n}v=\Res_{x}\frac{(1+x)^{\wt u+n}}{x^{2n+2}}Y(u,x)v
\end{eqnarray}
for $u,v\in V$ with $u$ homogeneous. The following was proved therein:

\bp{pdlmAn}
The subspace $O_{n}(V)$ is a two-sided ideal of 
the non-associative algebra $(V,*_{n})$ and the quotient algebra $V/O_{n}(V)$, denoted by $A_{n}(V)$,
is an associative algebra where ${\bf 1}+O_{n}(V)$
is an identity element and $\omega+O_{n}(V)$ is a central element denoted bar $\bar{\omega}$.
Furthermore, $A_{n}(V)$ admits an anti-automorphism $\phi$ given by $\phi(v)=e^{L(1)}(-1)^{L(0)}v$ for $v\in V$.
\ep

Let $(W,Y_{W})$ be a weak $V$-module and let $n\in \N$. Following \cite{dlm-anv}, set
\begin{eqnarray}
\Omega_{n}(W)=\{ w\in W\ |\  x^{n}Y_{W}(x^{L(0)}v,x)w\in W[[x]]\  \   \mbox{ for }v\in V\}.
\end{eqnarray}
It was proved in \cite{dlm-anv} that $\Omega_{n}(W)$ is an $A_{n}(V)$-module with $v$ acting as $v_{\wt v-1}$ for homogeneous $v\in V$. 
Suppose that $W=\oplus_{n\in \N}W(n)$ is an $\N$-graded weak $V$-module.
Then $\Omega_{n}(W)\supset \oplus_{i=0}^{n}W(i)$, which is a direct sum of $A_{n}(V)$-modules. 
Furthermore, if $W$ is irreducible with $W(0)\ne 0$, then for each $0\le i\le n$,
$W(i)$ is either zero or an irreducible $A_{n}(V)$-module, and $W(0),\dots, W(n)$ 
as $A_{n}(V)$-modules are non-isomorphic. 

The following was also obtained in \cite{dlm-anv} (Theorems 4.10 and 4.11):  

\bt{tdlm-anv}
A vertex operator algebra $V$ is rational if and only if $A_{n}(V)$ for all $n\in \N$ are finite-dimensional and semi-simple.
\et

It was also shown that $O_{n+1}(V)\subset O_{n}(V)$ for $n\in \N$ and the natural map 
$$\pi_{n}:\  A_{n+1}(V)\rightarrow A_{n}(V)$$
is an algebra epimorphism.

Here, we have:

\bt{tobservation}
Let $(V,\omega)$ be a vertex operator algebra. Suppose that $\omega'$ is a conformal vector of $V$ 
with integer conformal weights.
Then there exists a linear isomorphism $\Psi_{\omega,\omega'}:\ V\rightarrow V$ such that for each $n\in \N$,
$\Psi_{\omega,\omega'}(O_{n}(V,\omega))=O_{n}(V,\omega')$  and the induced linear map
\begin{eqnarray}
\Psi_{\omega,\omega'}^{(n)}: \   A_{n}(V,\omega)\rightarrow A_{n}(V,\omega')
\end{eqnarray}
is an algebra isomorphism.
Furthermore, if  $(V,\omega)$ is rational and if $(V,\omega')$ is a vertex operator algebra,
then $(V,\omega')$ is also rational. On the other hand, if $(V,\omega)$ is regular, 
then $(V,\omega')$ is a vertex operator algebra and it is regular.
\et

\begin{proof} For a general conformal vertex algebra $U$ of central charge $c$, 
Zhu (see \cite{zhu1}) proved that there is a new conformal vertex algebra structure on $U$ with
 the vertex operator map $Y[\cdot,x]$ defined by
$$Y[u,x]=Y(e^{xL(0)}u,e^{x}-1)\   \   \   \mbox{ for }u\in U$$
and with ${\bf 1}$ as the vacuum vector and
$\tilde{\omega}=\omega-\frac{1}{24}c{\bf 1}$ as the conformal vector. 
We denote this new conformal vertex algebra by $\exp(U,\omega)$. 
It was proved by Zhu and Huang (see \cite{zhu1}, \cite{huang}) that for any vertex operator algebra $(V,\omega)$,
$\exp(V,\omega)$ is isomorphic to $(V,\omega)$.

Let $n$ be any nonnegative integer. Recall from \cite{dlm-anv}  that for $u,v\in V$ with $u$ homogeneous, we have
\begin{eqnarray*}
&&u*_{n}v=\sum_{k=0}^{n}(-1)^{k}\binom{n+k}{n}\Res_{x}Y(u,x)v \frac{(1+x)^{\wt u+n}}{x^{n+k+1}},\nonumber\\
&&u\circ_{n}v=\Res_{x}\frac{(1+x)^{\wt u+n}}{x^{2n+2}}Y(u,x)v.
\end{eqnarray*}
By changing variable $x=e^{z}-1$ we get
\begin{eqnarray*}
u*_{n}v&=&\sum_{k=0}^{n}(-1)^{k}\binom{n+k}{n}\Res_{z}Y(u,e^{z}-1)v \frac{e^{z(\wt u+n+1)}}{(e^{z}-1)^{n+k+1}}\nonumber\\
&=&\sum_{k=0}^{n}(-1)^{k}\binom{n+k}{n}\Res_{z}Y(e^{zL(0)}u,e^{z}-1)v \frac{e^{z(n+1)}}{(e^{z}-1)^{n+k+1}}\nonumber\\
&=&\Res_{z}\left(\sum_{k=0}^{n}\binom{-n-1}{k}\frac{e^{z(n+1)}}{(e^{z}-1)^{n+k+1}}\right)Y(e^{zL(0)}u,e^{z}-1)v \nonumber\\
&=&\Res_{z}f_{n}(z)Y(e^{zL(0)}u,e^{z}-1)v
\end{eqnarray*}
and
\begin{eqnarray*}
u\circ_{n}v&=&\Res_{z}\frac{e^{z(n+1)}}{(e^z-1)^{2n+2}}Y(z^{zL(0)}u,e^{z}-1)v.
\end{eqnarray*}
On the other hand, from \cite{zhu1} we have $L[-1]=L(-1)+L(0)$.  
It then follows from Definition \ref{dtildeA} that 
\begin{eqnarray}
A_{n}(V,\omega)= \widetilde{A}_{n}(\exp(V,\omega)).
\end{eqnarray}
At the same time, this also proves that  $A_{n}(V,\omega')= \widetilde{A}_{n}(\exp(V,\omega'))$.
 As $\exp(V,\omega)$ and $\exp(V,\omega')$ are both isomorphic to $V$ as a vertex algebra, 
 there is a linear automorphism $\Psi_{\omega,\omega'}$ of $V$, which is a vertex algebra isomorphism 
 from $\exp(V,\omega)$ to $\exp(V,\omega')$.
Then $\Psi_{\omega,\omega'}$ gives rise to an algebra isomorphism 
 $$\Psi_{\omega,\omega'}^{(n)}:\  \widetilde{A}_{n}(\exp(V,\omega))\rightarrow \widetilde{A}_{n}(\exp(V,\omega'))$$ 
  for every $n\in \N$.
 Consequently, we have $A_{n}(V,\omega)\simeq A_{n}(V,\omega')$ for $n\in \N$. 
 If $(V,\omega')$ is a vertex operator algebra, it then follows immediately from
Theorem \ref{tdlm-anv} that $(V,\omega')$ is rational if  and only if $(V,\omega)$ is rational.

Suppose that $(V,\omega)$ is regular. From definition, each irreducible weak module is a $(V,\omega)$-module and there are finitely many irreducible weak modules up to isomorphism. By Theorem \ref{tdlm-anv}, $A_{n}(V,\omega)$ for $n\in \N$ 
are all finite-dimensional and semisimple. 
Let $W_1,\dots,W_{r}$ be a complete set of isomorphism-equivalence class representatives of irreducible weak $V$-modules, which are irreducible $(V,\omega)$-modules. 
From Zhu's  1-1 correspondence, $A(V,\omega)$ has exactly $r$ irreducible modules up to isomorphism, so does
$A(V,\omega')$ as $A(V,\omega')\simeq A(V,\omega)$.
This implies that there are $r$ irreducible $\N$-graded weak $(V,\omega')$-modules $M_1,\dots, M_r$ up to to equivalence. 
 These are $(V,\omega')$-modules since $\dim A_{n}(V,\omega')=\dim  A_{n}(V,\omega)<\infty$ for all $n\in \N$.
In particular, these are non-isomorphic irreducible weak $V$-modules. Thus
$M_1,\dots, M_r$ also form a complete set of isomorphism-equivalence class representatives of irreducible weak $V$-modules.
Consequently, $W_1,\dots,W_{r}$ are also $(V,\omega')$-modules.
Note that $V$ as a $(V,\omega)$-module is isomorphic to a finite direct sum of 
some $M_1,\dots, M_r$. It then follows that $V$ is also a $(V,\omega')$-module. Therefore, 
$(V,\omega')$ is a vertex operator algebra. The regularity assertion is clear.
 \end{proof}

\br{rnonvoa}
{\em Here, we show that it is possible that for a vertex operator algebra $V$, there exists a conformal vector
$\omega'$ with integer conformal weights such that $(V,\omega')$ is not a vertex operator algebra.
Let $\g$ be a finite-dimensional simple Lie algebra and let $\ell\in \C$
with $\ell \ne -h^{\vee}$ where $h^{\vee}$ is  the dual coxeter number of $\g$. We have a simple vertex operator algebra 
$L_{\hat{\g}}(\ell,0)$ where the conformal vector $\omega$ is obtained through the Segal-Sugawara construction.
Let $h$ be an ad-semisimple element of $\g$ with integer eigenvalues.
Set
\begin{eqnarray}
\omega_{h}=\omega-h(-2){\bf 1}=\omega-L(-1)h\in L_{\hat{\g}}(\ell,0)_{(2)}.
\end{eqnarray}
It is known that $\omega_{h}$ is a conformal vector with integer conformal weights, 
where 
$$L_{h}(n)=L(n)+(n+1)h(n)\   \    \   \mbox{ for }n\in \Z.$$
In particular, we have $L_{h}(-1)=L(-1)$ and $L_{h}(0)=L(0)+h(0)$.
Consider the special case with $\g=\sl(2,\C)$ and with 
$h$ taken to be the standard Cartan element. Set $\omega'=\omega_{h}$. 
The $L'(0)$-conformal weights of $e,f, h$ are
$$\wt' e=3,\   \   \   \wt' f=-1,  \   \   \   \wt' h=0.$$
Assume $\ell\notin \N$.
Note that the subalgebra of $\widehat{\sl_{2}}$ spanned by $f_{-1},e_{1},h(0)+{\bf k}$ is isomorphic to $\sl(2, \C)$ 
 and $e_{1}{\bf 1}=0,\  (h(0)+{\bf k}){\bf 1}=\ell {\bf 1}$. It follows that
  $(f_{-1})^{k}{\bf 1}\ne 0$ for all $k\in \N$.
 As $\wt' (f_{-1})^{k}{\bf 1}=-k$, we have that $L_{\hat{\g}}(\ell,0)_{(-k)}\ne 0$ for all $k\in \N$.
Thus $(L_{\hat{\g}}(\ell,0),\omega')$ is not  a vertex operator algebra.}
\er

Next, we discuss fusion rules.
Let $(V,\omega)$ be a conformal vertex algebra and let $W_1,W_2,W_3$ be weak $V$-modules. 
 An {\em intertwining operator of type} ${W_3\choose W_1W_2}$ (see  \cite{fhl}) is a linear map 
\begin{eqnarray*}
I(\cdot,x): &&W_1\rightarrow \left(\Hom (W_2,W_3)\right)\{x\}\\
&&w_1\mapsto I(w_1,x)=\sum_{\alpha \in \C}(w_1)_{\alpha }x^{-\alpha-1}\   \  (\mbox{where }(w_1)_{\alpha }\in \Hom (W_2,W_3)),
\end{eqnarray*}
 satisfying the conditions that for $w_1\in W_1,\ w_2\in W_2,\ \alpha\in \C$,  
 $$(w_1)_{\alpha+n}w_2=0\   \   \   \mbox{ for sufficiently large integers }n,$$
\begin{eqnarray}
&&x_0^{-1}\delta\left(\frac{x_1-x_2}{x_0}\right)Y(v,x_1)I(w_1,x_2)w_2-x_0^{-1}\delta\left(\frac{x_2-x_1}{-x_0}\right)I(w_1,x_2)Y(v,x_1)w_2
\nonumber\\
&&\hspace{3cm}=x_2^{-1}\delta\left(\frac{x_1-x_0}{x_2}\right)I(Y(v,x_0)w_1,x_2)w_2
\end{eqnarray}
for $v\in V,\ w_1\in W_1,\ w_2\in W_2$, and that
\begin{eqnarray}
I(L(-1)w_1,x)=\frac{d}{dx}I(w_1,x)\   \   \   \mbox{ for }w_1\in W_1.
\end{eqnarray}

Define a {\em lowest weight weak $V$-module} to be a weak $V$-module $W$ on which $L(0)$ acts semi-simply such that  
$W=\oplus_{n\in \N}W_{(\lambda +n)}$ for some $\lambda\in \C$  and $W_{(\lambda)}$ is an irreducible $A(V,\omega)$-module
which generates $W$ as a weak $V$-module. Note that if $V$ is of countable dimension, 
then from Lemma \ref {lzhu} every irreducible $\N$-graded weak $V$-module is a lowest weight weak $V$-module.

\bp{pfusion-rule}
Let $(V,\omega)$ be a vertex operator algebra and let $\omega'$ be any conformal vector of $V$, so that
$(V,\omega')$ is a conformal vertex algebra.
Assume that $W_2,W_3$ are weak $(V,\omega)$-modules and 
$W_1$ is a lowest weight weak $(V,\omega)$-module. Then there is a canonical 
isomorphism between the space of intertwining operators with respect to $\omega$ of type
${W_3\choose W_1W_2}$ and
 the space of intertwining operators with respect to $\omega'$ of the same type.
\ep

\begin{proof}  Set $a=\omega-\omega'\in V$. Then 
$$a_0=\omega_0-\omega'_{0}=L(-1)-L'(-1)=\D-\D=0\  \  \mbox{on }V.$$
In view of this, for any weak $V$-module $(W,Y_{W})$ we have 
$$[a_0,Y_{W}(v,x)]=Y_{W}(a_0v,x)=0\   \   \   \mbox{ for }v\in V.$$
Let $\lambda$ be the lowest $L(0)$-weight of $W_1$. Then 
$(W_1)_{(\lambda)}$ is an irreducible $A(V,\omega)$-module and $(W_1)_{(\lambda)}$ generates $W_1$.
As $[a_0,Y_{W_1}(\omega,x)]=0$, which implies $[a_0,L(0)]=0$,
we have $a_{0}(W_1)_{(\lambda)}\subset (W_1)_{(\lambda)}$.
Notice that $A(V,\omega)$ as a quotient space of $V$ is of countable dimension.
It then follows that $a_{0}$ acts on $(W_1)_{(\lambda)}$ as a scalar $\alpha\in \C$.
Since $(W_1)_{(\lambda)}$ generates $W_1$ and $[a_0,Y_{W_1}(v,x)]=0$ for all $v\in V$, 
it follows that  $a_0$ acts on $W_1$ as scalar $\alpha$.  That is, $L(-1)-L'(-1)=\alpha$ on $W_1$.

Now, suppose that $I(\cdot,x)$ is an intertwining operator of type ${W_3\choose W_1W_2}$ 
with respect to conformal vector $\omega$.
From definition, we have
$$I(L(-1)w,x)=\frac{d}{dx}I(w,x)\   \   \   \mbox{ for }w\in W_1.$$
Then 
\begin{eqnarray*}
e^{-\alpha x}I(L(-1)w,x)=e^{-\alpha x}\frac{d}{dx}I(w,x)=\frac{d}{dx}\left(e^{-\alpha x}I(w,x)\right)+\alpha e^{-\alpha x}I(w,x),
\end{eqnarray*}
which implies
\begin{eqnarray}
e^{-\alpha x}I(L'(-1)w,x)=\frac{d}{dx}\left(e^{-\alpha x}I(w,x)\right).
\end{eqnarray}
It can be readily seen  that $e^{-\alpha x}I(\cdot,x)$ also satisfies the Jacobi identity, so it is an intertwining operator 
with respect to conformal vector $\omega'$.
From this argument we also see  that if $I'(\cdot,x)$ is an intertwining operator with respect to conformal vector 
$\omega'$, then $e^{\alpha x}I'(\cdot,x)$ is an intertwining operator with respect to 
$\omega$. Now, the isomorphism assertion follows.
\end{proof}

\br{rhl}
{\em  We note that if $W_1,W_2,W_3$ are $(V,\omega)$-modules and $(V,\omega')$-modules as well, the assertion of 
Proposition \ref{pfusion-rule} also follows  from a theorem of \cite{hl} in terms of intertwining maps. }
\er

\end{document}